\documentclass[11pt]{amsart}

\usepackage[all]{xypic}

\usepackage{latexsym}

\usepackage{amssymb}
\usepackage{amsfonts}
\usepackage{amscd}
\usepackage{amsmath,amsthm}

\newtheorem{lemma}{Lemma}[section]

\newtheorem{theorem}[lemma]{Theorem}

\newtheorem{lem}[lemma]{Lemma}
\newtheorem{prop}[lemma]{Proposition}
\newtheorem{thm}[lemma]{Theorem}
\newtheorem{cor}[lemma]{Corollary}

{

}

\theoremstyle{definition}

\theoremstyle{remark}

\newtheorem{remark}[lemma]{Remark}

\numberwithin{equation}{section}

\newenvironment{pf}{\noindent{\bf Proof.}}{\hfill $\square$\medskip}

\def\AA{{\mathbb A}}
\def\CC{{\mathbb C}}

\def\NN{{\mathbb N}}

\def\PP{{\mathbb P}}

\def\RR{{\mathbb R}}

\def\ZZ{{\mathbb Z}}

\def\kol{{\bar k}}

\def\alphaol{{\bar{\alpha}}}

\def\tauol{{\bar \tau}}

\def\varphiol{{\bar \varphi}}

\def\sigmaol{{\bar \sigma}}
\def\betaol{{\bar{\beta}}}

\def\omegaol{{\bar \omega}}

\def\Rol{{\bar R}}

\def\0ol{{\bar 0}}
\def\1ol{{\bar 1}}
\def\2ol{{\bar 2}}
\def\ol2{{\bar 2}}
\def\3ol{{\bar 3}}
\def\4ol{{\bar 4}}
\def\5ol{{\bar 5}}
\def\6ol{{\bar 6}}
\def\7ol{{\bar 7}}
\def\8ol{{\bar 8}}
\def\9ol{{\bar 9}}

\def\bold0{{\bf 0}}
\def\bold1{{\bf 1}}
\def\bold2{{\bf 2}} 
\def\bold3{{\bf  3}}
\def\bold4{{\bf 4}}
\def\bold5{{\bf 5}}
\def\bold6{{\bf 6}}
\def\bold7{{\bf 7}}
\def\bold8{{\bf 8}}
\def\bold9{{\bf 9}}

\def\P2Skly{\PP^2_{Skly}}

\def\End{\operatorname {End}}
\def\Ext{\operatorname {Ext}}

\def\gr{\operatorname {gr}}

\def\Hom{\operatorname {Hom}}

\def\ker{\operatorname {ker}}

\def\th{\operatorname {th}}    
\def\Tor{\operatorname {Tor}}

\def\AAut{\operatorname{A\!\!\!Aut}}
\def\Aut{\operatorname{Aut}}

\def\dim{\operatorname{dim}}

\def\End{\operatorname{End}}

\def\ext{\operatorname{ext}}
\def\Ext{\operatorname{Ext}}

\def\Fract{\operatorname{Fract}}

\def\Gr{{\sf Gr}}

\def\hom{\operatorname{hom}}
\def\Hom{\operatorname{Hom}}

\def\id{\operatorname{id}}

\def\Iso{\operatorname{Iso}}

\def\Max{\operatorname{Max}}
\def\min{\operatorname{min}}

\def\Mod{{\sf Mod}}

\def\Pic{\operatorname{Pic}}

\def\rank{\operatorname{rank}}

\def\Spec{\operatorname{Spec}}

\def\sup{\operatorname{sup}}

\def\ul1{\operatorname{\underline{1}}}

\def\G{\mathop{\underline{\underline{\it \Gamma}}}\nolimits}

\def\l{\leftarrow}

\def\d{\downarrow}

\def\a{\alpha}
\def\b{\beta}
\def\c{\gamma}
\def\d{\delta}

\def\g{\gamma}
\def\l{\lambda}

\def\s{\sigma}

\def\ve{\varepsilon}

\def\D{\Delta}
\def\G{\Gamma}

\def\fa{{\mathfrak a}}

\def\fm{{\mathfrak m}}

\def\fp{{\mathfrak p}}

\def\sA{{\sf A}}

\def\sC{{\sf C}}

\def\sP{{\sf P}}

\def\wtf{{\widetilde{ f}}}

\def\wtm{{\widetilde{ m}}}

\def\wtalpha{{\widetilde{\alpha}}}

\def\whf{{\widehat{ f}}}
\def\whm{{\widehat{ m}}}

\def\cal{\mathcal}

\def\cD{{\cal D}}

\def\cF{{\cal F}}

\def\cL{{\cal L}}
\def\cM{{\cal M}}
\def\cN{{\cal N}}
\def\cO{{\cal O}}

\def\cS{{\cal S}}

\def\cX{{\cal X}}

\def\coh{{\sf coh}}

\def\Qcoh{{\sf Qcoh}}

\def\dirlim{\mathop{\vtop{\baselineskip -100pt\lineskip -1pt\lineskiplimit 0pt
\setbox0\hbox{lim}\copy0\hbox to \wd0{\rightarrowfill}}}\limits}
\def\invlim{\mathop{\vtop{\baselineskip -100pt\lineskip -1pt\lineskiplimit 0pt
\setbox0\hbox{lim}\copy0\hbox to \wd0{\leftarrowfill}}}\limits}

\def\I11{{1 \kern -0.8pt \! \mbox{l}}}
\def\mumu{{\mu\kern-4.2pt\mu}}
\def\bfmu{{\mu\kern-4.2pt\mu}}
\def\2slash{\backslash \! \backslash}

\def\boxtimes{\setbox0\hbox{$\Box$}\copy0\kern-\wd0\hbox{$\times$}}

\date{}                                           

\begin{document}

\title[The Weyl algebra and a quotient stack]
{A quotient stack related to the Weyl algebra}
\author{S. Paul Smith} 
\address{Department of Mathematics, Box 354350, University
of Washington, Seattle, WA 98195, USA}
\email{smith@math.washington.edu}

\thanks{The author was supported by NSF grant DMS-0602347}

\begin{abstract}
Let $A$ denote the ring of differential operators on the affine line with its two usual generators $t$
and ${{d}\over{dt}}$ given
degrees $+1$ and $-1$ respectively. Let $\cX$ be the stack having coarse moduli space the 
affine line $\Spec k[z]$ and isotropy groups $\ZZ/2$ at each integer point. Then the category of graded
$A$-modules is equivalent to the category of quasi-coherent sheaves on $\cX$.
\end{abstract}

\maketitle

\section{Introduction}

\subsection{}

Let $k$ be a field of characteristic zero. All vector spaces and algebras in this paper are taken over 
 $k$.

The first Weyl algebra is the ring $A=k\langle x,y\rangle/(xy-yx-1)$. We impose a $\ZZ$-grading 
on it by setting  $\deg x=1$ and $\deg y=-1$.
There is an isomorphism between $A$ and the ring
of differential operators with polynomial coefficients on the affine line $\Spec k[t]$  that is given by sending $x$ to ``multiplication by $t$'' and $y$ to $-d/dt$.

Our main result is that the category $\Gr A$ of $\ZZ$-graded $A$-modules  is equivalent to the 
category of quasi-coherent sheaves on a quotient stack $\cX$ whose coarse  moduli space is the affine line $\Spec k[z]$, and whose stacky structure consists of stacky points $B\ZZ_2$ supported at each integer point $n \in \ZZ \subset \AA^1_k$. 
We write
\begin{equation}
\label{equiv}
\Gr A \equiv \Qcoh\cX
\end{equation}
to denote this equivalence.

\subsection{}
We now  describe $\cX$.

Let $\ZZ_{\rm fin}$ be the group of finite subsets of $\ZZ$ with group operation given by ``exclusive or.''
Let $G$ be the affine group scheme whose coordinate ring is the group algebra $k\ZZ_{\rm fin}$
with its usual Hopf algebra structure. Since $\ZZ_{\rm fin}$ is $2$-torsion and is generated by the singleton sets $\{n\}$  a $k$-valued point $g\in G$ corresponds to a function $\ZZ \to \{ \pm 1\}$, $n \mapsto g(\{n\})$. 
We write $g_n$ for $g(\{n\})$. 

We define an action of $G$ on the ring
$$
C:=k\big[z\big]\big[\sqrt{z-n} \;\,  \big\vert \; n \in \ZZ\big]
$$
  by $g.\sqrt{z-n} := g_n \sqrt{z-n}$, and the  stack-theoretic quotient 
$$
\cX:=  \biggl[ {{\Spec C} \over{G}}\biggr].
$$
Its coarse moduli space is the affine line $\Spec k[z]$. 
Max Lieblich tells me that $\cX$ is an algebraic stack whose diagonal is
locally of finite type but not quasi-compact  (even though it is unramified).

\subsection{}

The action of $G$ on $C$ corresponds to the $\ZZ_{\rm fin}$-grading on $C$  given by
$$
\deg \sqrt{z-n}=\{n\}.
$$
A standard result for quotient stacks says that $\Qcoh \cX$ is equivalent to the category of $G$-equivariant
sheaves on $\Spec C$ or, equivalently, that there is an equivalence
$$
\Qcoh \cX \equiv \Gr(C,\ZZ_{\rm fin}),
$$
where $\Gr(C,\ZZ_{\rm fin})$ is the category of $\ZZ_{\rm fin}$ graded $C$-modules. 
Under this equivalence locally free coherent $\cO_\cX$-modules correspond to finitely generated projective
graded $C$-modules; for example,  $\cO_\cX$ corresponds to $C$.

We note that   $C$ is isomorphic to the polynomial ring $k[x_n \; | \; n \in \ZZ]$ 
modulo the relations $x_n^2+n=x_m^2+m$ for all $m$ and $n$, with the grading
given by $\deg x_n:=\{n\}$, and the isomorphism given by $z\leftrightarrow x_0^2$, and $\sqrt{z-n} \leftrightarrow x_n$. We therefore think of $k$-valued points of $\Spec C$ as elements in $k^\ZZ$.

The map $(a_i)_{i \in \ZZ} \mapsto (a_1 \ldots a_{2g+1}, a_0^2)$ is a surjective morphism from $\Spec C$
to a hyperelliptic curve of genus  $g$. If $k$ is not of characteristic two the fibers of this morphism are uncountable. When $k=\CC$ with its usual topology and $\CC^\ZZ$ is given the  product topology, 
and $\Max C$ is 
viewed as a subspace of $\CC^\ZZ$ with the subspace topology, the fibers are Cantor sets.

\subsection{}

At first sight, the equivalence (\ref{equiv}) is surprising. The Weyl algebra is an infinite dimensional 
$k$-algebra having no two-sided 
ideals other than zero and the ring itself so has no non-zero finite dimensional modules. 
As such it is ``very non-commutative.''  In stark contrast, $C$  is not only commutative but is even a graded PID meaning it is a domain and 
every graded ideal is  principal.  Graded right ideals in $A$ need not be principal.
Moreover, $C$ is a directed union of Dedekind domains. 
Although $C$ is not noetherian, it is noetherian from the graded 
perspective, meaning that $\Gr(C,\ZZ_{\rm fin})$ is a locally noetherian category. In particular, it has a set of noetherian generators.

\subsection{}
The relation between $A$, $C$, and $\cX$, can be viewed in the following way.
Let $\sA$ be the $k$-linear abelian category $\Gr  A$ but forget for the moment that it is $\Gr A$ and consider it
just as an abelian category. 
The endomorphism ring of the identity functor $\id_\sA$ is a 
polynomial ring in one variable. In the equivalences of $\sA$ with $\Gr A$, $\Gr C$, and $\Qcoh \cX$, 
$\End(\id_\sA)$ identifies with $A_0=k\big[t{{d}\over{dt}}\big]$, $C_\varnothing =k[z]=k[x_0^2]$, and $\G(\cX,\cO_\cX)$, respectively.
The Picard group, $\Pic\sA$, of $\sA$ is defined as the group of auto-equivalences modulo isomorphism. 
(It {\it does} contain the usual Picard group of $\cX$ where one identifies an invertible $\cO_\cX$-module $\cL$ with the auto-equivalence $\cL \otimes -$.)
Thinking of elements of $\Pic\sA$ as being like invertible sheaves, or line bundles, one is led to 
associate to each subgroup $\G$ of  $\Pic \sA$ a ``homogeneous coordinate ring''
$$
\bigoplus_{F \in \G}\Hom(\id_\sA,F)
$$
where $\Hom(\id_\sA,F)$ consists of natural transformations from $\id_\sA$ to $F$. This point of view lies at the heart of the work of Artin and Zhang's conception of non-commutative algebraic geometry \cite{AZ}.
There are particular subgroups of $\Pic\sA$ isomorphic to $\ZZ_{\rm fin}$ and $\ZZ$, and the 
``homogeneous coordinate rings'' associated to these subgroups are 
isomorphic to $C$ and $A$, respectively. The subgroup isomorphic to $\ZZ_{\rm fin}$ also 
identifies with $\Pic\cX$. In some sense, $x_n$, or rather the result of its action on
$\id_\sA$, is an endo-functor of $\sA$ that is a square root of the endo-functor of $\Gr A$ that is given by a 
left action  of the operator $t{{d}\over{dt}}-n$ on graded right $A$-modules (see sections \ref{sect.preps} and 
\ref{sect.R0.action}).  

The subgroup of $\Pic(\Gr A)$ isomorphic to $\ZZ_{\rm fin}$ was found by Sue Sierra \cite{Sue1}.

\subsection{}
The stimulus for this paper was Sierra's work on the graded Weyl algebra \cite{Sue1} 
and especially her  ``picture''
\begin{equation}
\label{Sues.pic}
\UseComputerModernTips
  \xymatrix{
\cdots &\ar@{-}[l] \!\! \colon \! \! \! \!  \ar@{-}[r] &  \!\! \colon \! \! \! \!  \ar@{-}[r] & \!\! \colon \! \! \! \!  \ar@{-}[r] & \!\! \colon \! \! \! \!  \ar@{-}[r]  &\!\! \colon \! \! \! \! \ar@{-}[r]& \cdots
 }
 \end{equation}
  of the simple graded $A$-modules which reminded this author of a stack on the affine line with stacky structure $B\ZZ_2$ at each integer point.
  Each point  in Sierra's picture represents a simple graded $A$-module: if $\l \in k-\ZZ$ there is a 
  single simple graded  $A$-module up to isomophism, namely $A/(xy-\l)A$;
   if $n \in \ZZ$ there are two simple modules, 
   $$
   X(n):=\biggl({{A}\over{xA}} \biggr)(n) 
   \quad \hbox{and} \quad   
   Y(n):=\biggl({{A}\over{yA}} \biggr)(n-1).
   $$
   (The isomorphism between $A$ and $\cD(\AA^1)$ may be chosen 
    so that $Y(1)$ 
 corresponds to the natural module $k[t]$ and $X(1)$  corresponds to the module
$k[t,t^{-1}]/k[t]$.)
    There is a non-split extension of  $X(n)$ by $Y(n)$ and a non-split extension of $Y(n)$ by $X(n)$ for each $n$.

    The underlying line in (\ref{Sues.pic})  should be thought of as $\Spec k[xy]$ and the two points at $n \in \ZZ$
    represent, in some sense, the two formal square roots of $-xy - n$ (because 
    $-xy$ corresponds to $t{{d}\over{dt}}$).

    We call the two points at $n \in \ZZ \subset k$ ``fractional points''. There are two reasons for this.     
 First,  if $n \in \ZZ$ and $\l \in k-\ZZ$ there is an equality $[X(n)]+[Y(n)]= [A/(xy-\l)A]$  in the Grothendieck group of $\Gr A$, and under the equivalence with  $\Qcoh \cX$, $ [A/(xy-\l)A]$  identifies with the skyscraper 
 sheaf $\cO_\l$ at the point $\l \in \Spec k[z]$.
  Second, there is some 
 consistency with the  notion of a ``fractional brane'' or ``brane fractionation'', where a brane represented by a 
 point in the Azumaya locus ``fractionates'' when it moves to the non-Azumaya locus.  
 
    Sierra's picture can also be viewed as a depiction of the stack 
  $\cX$. The line given by collapsing the ``fractional points'' is the coarse moduli space  $\Spec k[z]$ 
  of $\cX$ and the two points at $n$ correspond to the skyscraper sheaf $\cO_n=k[z]/(z-n)$ 
  endowed with the  trivial and sign representations of the isotropy 
  group at $n$, and $\l \in k-\ZZ$ corresponds to $\cO_\l=k[z]/(z-\l)$.  
If $\chi_{sgn}$ and $\chi_{triv}$ denote the sign and trivial representations of  the appropriate isotropy groups, then under the equivalence of categories $\Gr A \equiv \Qcoh \cX$ there are correspondences
  \begin{align*}
 X(n)  \leftrightsquigarrow  &
		 \begin{cases}
			  \cO_n \otimes   \chi_{triv} & \text{  if $n\le 0$}
			 \\
			 \cO_n \otimes   \chi_{sgn} & \text{  if $n \ge 1$}
		\end{cases}
 \\
  Y(n)  \leftrightsquigarrow  & \; 
  		 \begin{cases}
			  \cO_n \otimes   \chi_{sgn} & \text{  if $n\le 0$}
			 \\
			 \cO_n \otimes   \chi_{triv} & \text{  if $n \ge 1.$}
		\end{cases}
 \end{align*}
 Under the direct image functor  for the morphism from $\cX$ to its coarse moduli space ``half'' the $X(n)$s
 and  ``half'' the $Y(n)$s are sent to zero.

Under the equivalence with $\Gr C$, $X(n)$ corresponds to $C/(x_n)$ when $n \le 0$ and to 
$Cx_n/(x_n^2)$ when $n \ge 1$; similarly, $Y(n)$ corresponds to $C/(x_n)$ when $n \ge 1$ and to 
$Cx_n/(x_n^2)$ when $n \le 0$. 
\subsection{}

 Further evidence of a possible relation between $\Gr A$ and $\Qcoh \cX$ is the behavior of the Ext-groups between simple modules. Sierra showed that  the only non-trivial
 extensions between non-isomorphic simple graded $A$-modules are the following: 
 For all $i,j \in \ZZ$
 \begin{equation}
 \label{eq.ext.gps}
 \ext^1_A(X(i),Y(j)) \cong \ext^1_A(Y(j),X(i)) \cong 
 	\begin{cases}
	k & \text{if $i=j$, and}
	\\
	0 & \text{if $i \ne j$}
	\end{cases}
\end{equation}
where $\ext^1_A$ denotes $\Ext^1_A$  in $\Gr A$ \cite[Lem. 4.3]{Sue1}. This is ``the same'' as the behavior of the ext-groups between the simple objects in $\Qcoh \cX$.

\subsection{}

Using the equivalence between $\Gr A$ and $\Qcoh \cX$, 
the direct image functor for the morphism  from $\cX$ to its coarse moduli space transfers to a 
functor $\Gr A \to \Mod k[z]$. That functor sends a graded $A$-module to its degree zero component.
For example, $A$ viewed as a graded right $A$-module is sent to $k[xy]$ which identifies with 
$k[z]$. We therefore write $z$ for the element $xy$ of $A$ and think of it both as an element of $A$ and 
as the coordinate function on  the affine line that is the coarse moduli  space for $\cX$, i.e., 
$k[z] =C^G =C_{\varnothing}$. 
We note that  $k[z]$  is equal to $k[x_n^2]$ for all $n \in \ZZ$.

  \subsection{}

Having obtained the equivalence between $\Gr A$ and $\Qcoh \cX$ or, equivalently, with
$\Gr(C,\ZZ_{\rm fin})$, one can obtain alternative proofs of many of Sierra's results by transferring 
results from $\Gr(C,\ZZ_{\rm fin})$ to $\Gr A$. 
This is a good thing because the fact that $C$ is a graded PID makes the study of its graded modules quite
straightforward. 

To illustrate this point we compute the Grothendieck and Picard groups of  $\Gr(C,\ZZ_{\rm fin})$
directly using $C$ rather than using the equivalence of categories and quoting Sierra's result that computes
those invariants for $\Gr A$.

  \subsection{}
Section \ref{sect.prelim} is preparatory, setting up notation, and recalling some well-known facts. 
Section \ref{sect.autom.twist} concerns the Picard group $\Pic(R,\G)$  of the category 
$\Gr(R,\G)$ of graded modules over a ring $R$ graded by an abelian group $\G$. The results there
 may be of independent interest. The notion of an ``almost-automorphism'' of $(R,\G)$ is introduced and we 
show that every almost-automorphism  determines an autoequivalence of $\Gr(R,\G)$. 
Proposition \ref{prop.wisom} shows that  a pair of 
autoequivalences $F$ and $G$ such that $F(R(i)) \cong G(R(i))$ for all $i \in \G$
are naturally isomorphic if the endomorphism ring of every homogeneous component $R_i$  
is isomorphic to $R_0$. There is a group homomorphism $\Pic(R,\G) \to \Aut \End(\id_{\Gr(R,\G)})$.
In Proposition \ref{prop.R0} a criterion, which is satisfied by $(A,\ZZ)$ and $(C,\ZZ_{\rm fin})$,
 is given that implies there is a ring isomorphism $R_0 \to \End(\id_{\Gr(R,\G)})$. 
In particular, in this situation every
graded right $R$-module can be given the structure of an $R_0$-$R$-bimodule. 

Section \ref{sect.C.ring} concerns the structure of  $C$ as an ungraded ring and also examines the 
``variety'' $X \subset \CC^\ZZ$ of which $C$ is the coordinate ring. By 
definition, $X$ is the zero locus of the equations $x_n^2+n = x_m^2+m$, $m,n \in \ZZ$. The 
topological structure of $X$ is examined when $\CC^\ZZ$ is given various topologies and $X$ is
given the subspace topology. For example, when $\CC^\ZZ$ is given the product topology the fibers
of each projection $x_n:X \to \CC$ are Cantor sets. With the box topology $X$ becomes discrete.
With an appropriate embedding in $\ell^\infty(\ZZ)$, $X$ has uncountably many 
connected components, all homeomorphic to one another and permuted by the 
action of the group $\{\pm 1\}^\ZZ$ acting by coordinate-wise multiplication. 
Each component can be given the structure of a Riemann
surface with respect to which the coordinate functions $x_n$ are holomorphic.

Section \ref{sect.C.gr} establishes the properties of $C$ as a graded ring and culminates in 
a proof that the categories $\Gr(A,\ZZ)$ and $\Gr(C,\ZZ_{\rm fin})$ are equivalent. 
However, that equivalence is {\it not} used in sections \ref{sect.simples}---\ref{sect.symm}. 
Thus, the results in those sections are independent of Sierra's work, and provide alternative
proofs of several of her results. 

Section \ref{sect.simples} classifies the simple graded $C$-modules and focuses on those 
that are supported at the stacky points on $\cX$. We call those {\it special}. They correspond to the graded $A$-modules labelled $X(n)$ and $Y(n)$ above. As in Sierra's analysis, 
they play a central role in this paper. For example,  Corollary \ref{cor.simples.funct} shows that an autoequivalence 
of $\Gr (C,\ZZ_{\rm fin})$ is determined up to isomorphism by its action on the isomorphism classes 
of the special simple modules. 

The special simples may be characterized as those simple graded modules $S$
such that $\ext^1(S,S')$ is non-zero for some simple module $S'$ that is not isomorphic to $S$.
They may also be characterized as those simple graded modules $S$ such that $\hom(P,S)=0$ 
for some non-zero projective graded module (Proposition \ref{prop.specials}). 
Sierra exploits the first characterization in her paper
whereas we choose to exploit the second characterization in this paper so as to provide a
different perspective. 

Section \ref{sect.K0} computes the Grothendieck group of the category of finitely generated graded $C$-modules or, equivalently, that of the $G$-equivariant locally free sheaves on $\Spec C$. We compute
that Grothendieck group as an explicit quotient ring of the group algebra $\ZZ\ZZ_{\rm fin}$. In passing, we 
prove that the isomorphism classes of finitely generated projective graded $C$-modules or, equivalently,  locally free $G$-equivariant sheaves on $\Spec C$, are in natural bijection with the 
finite multi-sets of integers.

Section \ref{sect.symm} shows that the translation and reflection symmetries of Sierra's picture (\ref{Sues.pic})
can be implemented at the functorial level by the functor $\tau_*$ induced by the automorphism 
$\tau:z \mapsto z+1$, or $x_n \mapsto x_{n-1}$, of $C$, and by the functor $\varphi_*$ induced by the 
almost-automorphism $\varphi:x_n \mapsto x_{-n}, \;
x_n^2 \mapsto -x_{-n}^2, \; z \mapsto -z$,  respectively. The reflection symmetry cannot be induced at the functorial level by an automorphism of $C$ unless $k$ contains $\sqrt{-1}$. The main result in  section
 \ref{sect.symm}  is the computation of $\Pic(C,\ZZ_{\rm fin})$. We show it fits into a sequence
$$
1 \to \ZZ_{\rm fin} \to \Pic(C,\ZZ_{\rm fin}) \to \Iso(\ZZ) \to 1
$$
where $\Iso(\ZZ)$ is the isometry group of $\ZZ$, the infinite dihedral group, abstractly.

Section \ref{sect.A+C} makes a direct comparison between $\Gr (A,\ZZ)$ and $\Gr (C,\ZZ_{\rm fin})$. 
Proposition \ref{prop.match.simples} shows how the special simples over each ring correspond
under the equivalence of categories---the correspondence is not what the notation might lead one to 
expect. Theorem \ref{thm.iotaJ.CJ} shows that  the auto-equivalences $\iota_J$, $J \in \ZZ_{\rm fin}$, 
found by Sierra correspond  to the Serre twists $(J)$ on $\Gr (C,\ZZ_{\rm fin})$. 
Proposition \ref{prop.shift.Sigma} shows that the Serre twist, $(1)$, on $\Gr (A,\ZZ)$  
corresponds to the auto-equivalence $(\{1\}) \circ \tau_*$  of $\Gr (C,\ZZ_{\rm fin}) $. 
Proposition \ref{prop.fourier} shows that the auto-equivalence of $\Gr (A,\ZZ)$ induced by the automorphism $x \mapsto y$ and $y \mapsto -x$ corresponds to the 
auto-equivalence $\tau_*\varphi_*$ of $\Gr (C,\ZZ_{\rm fin})$  induced by the almost-automorphism $\tau\varphi$.  

The equivalence between $\Gr (A,\ZZ)$ and $\Gr (C,\ZZ_{\rm fin})$ was proved in section \ref{sect.C.gr}
 by starting with $C$ and then showing that $A$ was the endomorphism ring of a certain
 bigraded $P$-module. In section \ref{sect.tw.hcr} we take the opposite approach and 
 show that $C$ can be constructed from  $\Gr (A,\ZZ)$ as a sort of twisted homogeneous coordinate ring.

\subsection{}
The results about graded $C$-modules can be translated into results about $\Qcoh \cX$ or, equivalently, about the $G$-equivariant sheaves on $\Spec C$. For example 
\begin{enumerate}
  \item 
  the Grothendieck group $K_0(\cX)$ is a free abelian group of countable rank, and we 
  present it as an explicit quotient of the group algebra $\ZZ\ZZ_{\rm fin}$ in Theorem \ref{thm.K0.relns};
  \item{}
  the invertible $\cO_\cX$-modules are, up to isomorphisms, the twists of $\cO_\cX$ by the characters of 
  $G$ or, equivalently,  $\Pic \cX \cong \ZZ_{\rm fin}$ (Corollary \ref{cor.PicX});
  \item 
   every locally free $\cO_{\cX}$-module is a direct sum of invertible $\cO_{\cX}$-modules
   (Proposition \ref{prop.pid.mods});
  \item{}
   the locally free $\cO_\cX$-modules are, up to isomorphisms, in natural bijection with 
   the finite multi-sets of integers.  
\end{enumerate}

\subsection{Acknowledgements}

I am very grateful to Sue Sierra for answering my questions about her work.
I am grateful to Robin Graham for  telling me the results  in Propositions \ref{prop.box} and \ref{prop.RS}, and to Lee Stout for telling me the result in Proposition \ref{prop.not.RS}. 
I also thank Jack  Lee and John Sylvester for useful discussions.

  \section{Preliminaries}
  \label{sect.prelim}

\subsection{The Weyl algebra $A$}

The $\ZZ$-grading on $A$ given by
$$
\deg x=1 \quad \deg y=-1
$$
is sometimes called the {\sf Euler grading} because $A_n$ consists of those operators/elements 
$a$ such that
$[D,a]=na$ where $D$ is the Euler vector field/derivation
$$
D=t{{d}\over{dt}}.
$$

\subsection{The twist functor on $\Gr A$}

For each $n \in \ZZ$, we define Serre's twist automorphism $M \mapsto M(n)$ on $\Gr A$ by declaring 
that $M(n)$ is equal to $M$ as a right $A$-module but the grading is now given by
$$
M(n)_i = M_{n+i}.
$$
{\it This notation differs from Sierra's: }
her primary twist functor is denoted by $M \mapsto M\langle n\rangle$ where $M\langle n\rangle=M$
and
$$
M\langle n\rangle_i = M_{i-n}.
$$ 
Thus $\langle n\rangle =(-n)$.

\subsection{The $\ZZ_{\rm fin}$-grading on $C$}

Let $\ZZ_{\rm fin}$ be the group of finite subsets of the integers with group operation
$$
I \oplus J :=I \cup J - I\cap J  =(I-J) \cup (J-I).
$$
The identity is the  empty set $\varnothing$.  It is easy to see that  $\ZZ_{\rm fin}$ is  the direct sum of the 
two-element subgroups  $\{\varnothing,\{i\}\}$, $i \in \ZZ$. 

Define the commutative ring 
\begin{align*}
C:= & k[ x_n \; | \; n \in \ZZ]   \qquad \hbox{modulo the relations}
\\
x_n^2 + n  = & x_m^2+m, \qquad  \hbox{for all } n,m \in \ZZ.
\end{align*}
We make $C$ a $\ZZ_{\rm fin}$-graded ring by declaring that
$$
\deg x_n=\{n\}.
$$
The identity component, $C_\varnothing$, is equal to $k[x_n^2]$ for all $n \in \ZZ$.
 For each $I \in \ZZ_{\rm fin}$ we define
 $$
 x_I:=\prod_{i \in I} x_i
 $$
 with the convention that $x_\varnothing=1$.\footnote{Whenever things are indexed by 
 elements of $\ZZ_{\rm fin}$ we write $x_i$ rather than $x_{\{i\}}$ for the element indexed by  $\{i\}$. } 
The homogeneous components of $C$ are $C_I:=C_\varnothing  x_I$. We will see below that $C$ is a domain
so each $C_I$ is isomorphic to $C_\varnothing $ as a $C_\varnothing $-module.

\subsubsection{The ring $\ZZ_{\rm sub}$} 

Let $\ZZ_{\rm sub}$ denote the set of all subsets of $\ZZ$.
Then $\ZZ_{\rm sub}$ is a commutative ring with identity with respect to the product given by 
intersection and addition given by $\oplus$. Its identity is $\ZZ$ and its zero element is $\varnothing$. 
Thus $\ZZ_{\rm fin}$ is a subring of $\ZZ_{\rm sub}$ without identity.  Since every element of $\ZZ_{\rm sub}$ 
is idempotent, $\ZZ_{\rm sub}$ is a Boolean ring.

We note the identity $I-J=I\oplus (I \cap J) = I \cap (\ZZ \oplus J)$.

 \subsubsection{Notation for $\ZZ_{\rm fin}$} 
 If $J \in \ZZ_{\rm fin}$ and $n \in \ZZ$ we adopt the notation:
 \begin{itemize}
  \item 
  $
 n+J:=\{n+j \; | \; j \in J\};
 $
  \item 
   $
 nJ:=\{nj \; | \; j \in J\}.
 $
\end{itemize}

\subsection{Categories of graded modules}
\label{sect.gr.rings}

If $R$ is a ring graded by a group $\G$ we write $\Gr(R,\G)$ for the category of $\G$-graded  right 
$R$-modules with degree-preserving homomorphisms.  If $\G=\ZZ$ we often write $\Gr R$ for 
$\Gr(R,\ZZ)$.

We will write $\hom_R(M,N)$ for the degree preserving $R$-module homomorphisms from one 
 $\G$-graded $R$-module $M$ to another $N$. We denote the right derived functors of 
 $\hom_R$  by $\ext^i_R$, $i \ge 0$. 

There is a category of graded rings in which the objects are pairs $(R,\G)$ consisting of a group $\G$
and a $\G$-graded ring $R$. Morphisms are pairs $(\a,\alphaol):(R,\G) \to (R',\G')$ where $\alphaol:
\G \to \G'$ is a group homomorphism and $\a:R \to R'$ is a ring homomorphism such that $\a(r_i) \subset R'_{\alphaol i}$ for all $i \in \G$.

Associated to $(\a,\alphaol)$ is a functor $\a^*:\Gr(R,\G) \to \Gr(R',\G')$ and its right adjoint 
$\a_*:\Gr(R',\G') \to \Gr(R,\G)$.  
If $M \in \Gr(R',\G')$, then $(\a_* M)_i:=M_{\alphaol i}$ for all $i \in \G$ and $R$ acts on $\a_*M$   
via the homomorphism $\a$.

If $m \in M_{\alphaol i}$ we will  label it as $\a_*m$ when we think of it as an element in $\a_*M$.
Hence $x.\a_*m=\a_*(\a(x) m)$.

\subsection{The affine group scheme $G$}
\label{sect.gp.sch}

The  group algebra of $\ZZ_{\rm fin}$ is
$$
k\ZZ_{\rm fin} = {{k[u_i \; | \; i \in \ZZ]}\over{(u_i^2-1 \; | \; i \in \ZZ)}}
$$
where $u_i$ is an alias for the element $\{i\}$. The group algebra is given its usual
Hopf algebra structure and we define  the  affine group scheme
$$
G:=\Spec k\ZZ_{\rm fin}.
$$ 
The letter $G$ will denote both the group scheme and the group of $k$-valued points of 
 $\Spec k\ZZ_{\rm fin}$. 
 
 Let $\{\pm 1\}^{\ZZ}$ be the group of functions $\ZZ \to \{\pm 1\}$. 
 There is an isomorphism $G\to \{\pm 1\}^{\ZZ}$ that sends $g \in G$ 
to the function $\ZZ \to \{\pm 1\}$ given by $i \mapsto u_i(g)$. 
 Thus, $G$ is isomorphic to a countable product of copies of 
$\{\pm 1\}$. When we consider an element $g \in G$ as a function $\ZZ \to \{\pm 1\}$
that function will always be given by $g(i)=u_i(g)$.

\subsection{The stack}

We define an algebraic action of $G$ on $\Spec C$ by declaring that $g \in G$ acts on $x_i$ by 
$$
g.x_i=g(i) x_i,
$$
and define the quotient stack
$$
\cX:=  \biggl[ {{\Spec C} \over{G}}\biggr].
$$
The category of quasicoherent $\cO_{\cX}$-modules is denoted by $\Qcoh \cX$ and 
 is equivalent to the category of $G$-equivariant $C$-modules---we will usually think of it in this way. 

The invariant subring of $C$ is $C^G=C_\varnothing =k[x_n^2]$ for all $n$. The coarse moduli space of $\cX$ is therefore the affine line $\Spec k[x_0^2]$.

We will denote $k$-valued points in $\Spec C$ by tuples $(a_i) \in k^{\ZZ}$, $a_i \in k$, such a point corresponding to the maximal ideal 
$$
\sum_{i \in \ZZ}  (x_i-a_i).
$$ 
The relations for $C$ imply that  at most one $a_i$
is zero. Therefore the points  having a non-trivial isotropy group are those for
which one $a_i$ is zero.  The isotropy group at such a point is isomorphic to $\ZZ_2$. Such points
are those where $x_0^2$ takes an integer value. Hence all the stacky structure on $\cX$
occurs at the integer points $x_0^2=n$, $n \in \ZZ$, on the coarse moduli space $\Spec k[x_0^2]$.

\subsection{} 

Because  $\cO(G)$ is the group algebra 
$k\ZZ_{\rm fin}$ a rational representation of $G$ is the
 same thing as a $\ZZ_{\rm fin}$-graded  vector space. 
 In particular, the $\ZZ_{\rm fin}$-grading on $C$ is that induced by the action of $G$ on $\cO(C)$.  

It is a standard result that the category of $G$-equivariant $C$-modules is equivalent to the category of 
$\ZZ_{\rm fin}$-graded $C$-modules, so there is an equivalence of categories 
$$
\Qcoh \cX \equiv \Gr(C,\ZZ_{\rm fin})
$$
where the latter denotes the category of $\ZZ_{\rm fin}$-graded $C$-modules.

Our main result, namely that 
$$
\Gr A \equiv \Qcoh \cX ,
$$
 will be proved by showing that $\Gr A \equiv \Gr(C,\ZZ_{\rm fin})$.

\section{Autoequivalences of categories of graded modules}
\label{sect.autom.twist}

The {\sf Picard group} of a category is its group of auto-equivalences  modulo
natural isomorphism. 

Let $(R,\G)$ be a graded ring.
We write $\Pic(R,\G)$ for the Picard group of $\Gr(R,\G)$. 
As remarked in section \ref{sect.gr.rings}, an automorphism $(\a,\alphaol)$ of $(R,\G)$ induces an automorphism $\a_*$  of $\Gr(R,\G)$.  We write $[\a_*]$ for the isomorphism class of $\a$.
This passage $\a  \rightsquigarrow [\a_*]$ can be mimicked for
certain maps $\a:R \to R$ that are not automorphisms.  

\subsection{Almost-automorphisms of $(R,\G)$}
\label{sect.almost-autom}

Let $k$ be a   subring of $R_0$ that is central in $R$. 
Let $k^\times$ denote the group of units in $k$.

An {\sf almost-automorphism} of $(R,\G)$ is a triple $(\a,\alphaol,\l)$ consisting of 
\begin{enumerate}
  \item 
 an automorphism    $\alphaol$ of $\G$,
  \item 
 a  $k$-module automorphism   $\a:R \to R$ of $R$ 
 such that  $\a(R_i) = R_{\alphaol i}$ for all $i$, and 
  \item 
  a normalized 2-cocycle $\l:\G \times \G \to k^\times$, $(i,j) \mapsto \l_{ij}$,   i.e., 
  $
 \l_{00}=1
 $
 and 
\begin{equation}
\label{eq.aa}
\l_{h,i+j}\l_{ij}=\l_{h+i,j}\l_{hi}
\end{equation}
for all $h,i,j \in \G$, such that 
$$
\a(xy)=\l_{hi} \a(x) \a(y)
$$
for all $x \in R_h$ and all $y \in R_i$ and all $h,i \in \G$. 
\end{enumerate}

\subsubsection{}
Because $\l_{00}=1$, the restriction of $\a$ to 
$R_0$ is a $k$-algebra automorphism of $R_0$.

\subsubsection{}

An automorphism $(\a,\alphaol)$ of $(R,\G)$ is an almost-automorphism with $\l_{hi}=1$ for all $h$ and $i$.

\subsubsection{}
\label{sect.beattie}
I am grateful to Margaret Beattie for the following observation.

Suppose $\l$ is a normalized 2-cocycle and $(\a,\alphaol)$ a pair satisfying conditions (1) and (2) in section \ref{sect.almost-autom}. Let $\nu=\l \circ (\alphaol^{-1} \times \alphaol^{-1})$. Then $\nu$ is also a normalized 2-cocycle and there is a standard construction of a new graded $k$-algebra $(R^\nu,\G)$ which is $(R,\G)$ as a graded $k$-module, but endowed with a new multiplication
$$
x*y:=\nu_{ij} xy
$$
for $x \in R_i$ and $y \in R_j$. Beattie observed that $(\a,\alphaol,\l)$ is an almost 
automorphism of $R$ if and only if
$(\a,\alphaol)$ is an isomorphism of graded $k$-algebras $(R,\G) \to (R^\nu,\G)$.

\begin{lem}
The set of $k$-linear almost-automorphisms of a graded $k$-algebra $(R,\G)$
form a group with respect to the product
$$
(\a,\l) * (\b,\nu):=(\a\b,\xi),  \quad \hbox{ where $ \;  \xi_{ij}:= \l_{\betaol i,\betaol j} \nu_{ij}$ 
for all $i$ and $j$. }
$$
\end{lem}
\begin{pf}
It is straightforward to check that $(\a\b,\xi)$ is an almost-automorphism. The identity
automorphism $(\id_R,\id_\G)$ has the property that $(\id_R,1)*(\a,\l) =(\a,\l)*(\id_R,1)=(\a,\l)$,
so is an identity for the product $*$.

If $(\a,\l)$ is an almost-automorphism so is $(\a^{-1},\zeta)$ where
$$
\zeta_{ij}:=\l_{\alphaol^{-1}i,\alphaol^{-1}j}^{-1}
$$
for all $i$ and $j$. It is easy to see that $(\a,\l) * (\a^{-1},\zeta)$  and
$(\a^{-1},\zeta)  * (\a,\l)$ are equal to  $(\id_R,1)$. 
The set of almost-automorphisms is therefore a group. 
\end{pf}

We write $\AAut(R,\G)$ for the group of almost automorphisms of $(R,\G)$.

\subsection{}

For each $(i,j) \in \G \times \G$, let $\Rol_{i,j}:=R_{i-j}$ and define
$$
\Rol := \bigoplus_{i,j \in \G} \Rol_{i,j}.
$$
The components $\Rol_{i,j}$ can be viewed as the sets of morphisms $j \to i$ in the category $\sC(R,\G)$ whose objects are the elements of $\G$, and in which composition of 
morphisms is given by multiplication in $R$. Thus $\Rol$ is a ring but does not have an identity if $\G$ is infinite, though it does have ``local'' units.

Let $(\a,\alphaol,\l)$ be an almost-automorphism of $(R,\G)$. 
Define a $k$-linear map $\wtalpha:\Rol \to \Rol$ by 
$$
\wtalpha(x) := \l_{i-j,j}\a(x)
$$
for $x \in \Rol_{i,j}$. 
Then  $\wtalpha$ is an algebra  automorphism of $\Rol$.

We define an automorphism $F=F_{\a,\alphaol,\l}:\sC(R,\G) \to \sC(R,\G)$  of $\sC(R,\G)$ by declaring that 
$Fi: = \alphaol i$, $i \in \G$, and on a morphism $x \in \Rol_{i,j}$  its action is
$$
x \mapsto Fx:= \l_{i-j,j}\a(x).
$$

\subsection{Automorphisms of $\Gr(R,\G)$ induced by almost-automorphisms.}
\label{sect.aauts+auts}

Let $\a$ be an almost-automorphism of $(R,\G)$.
 Let $M \in \Gr(R,\G)$. We define $\a_*M$ to be $M$ endowed with the grading
 $$
 (\a_* M)_h := M_{\alphaol h}.
 $$
 We write $\a_*m$ for an element  $m\in M$ viewed as an
 element of $\a_* M$. It is easy to see that $\a_*M$ becomes a graded $R$-module when the action of 
 $x \in R_i$ on an element $\a_*m \in (\a_*M)_h$ is defined to be
 $$
 (\a_* m).x :=\a_*(m\a(x))  \l_{hi} .
 $$
 The only point to be checked is the associative law, $((\a_*m).x).y= (\a_*m).(xy)$, which follows from the 
 identity (\ref{eq.aa}).
 
 \begin{lem}
 \label{lem.almost-autom}
 Let $\a$ be an almost-automorphism of $(R,\G)$. If $f:M \to N$ is a homomorphism of graded $R$-modules, then the map $\a_*f$ defined by 
 $$
 (\a_* f)(\a_* m):= \a_*(fm)
 $$
 is a homomorphism $\a_*M \to \a_* N$. With these definitions, $\a_*$ becomes an automorphism of 
 $\Gr(R,\G)$.
 \end{lem} 
 
 \subsubsection{}
 \label{sect.beattie.2}
 There is another way to view the auto-equivalence $\a_*$ associated to an almost automorphism by 
 using Beattie's observation in section \ref{sect.beattie}. First, there is an
  identification $\Gr(R,\G) = \Gr(R^\nu,\G)$ because every left $R$-module $M$ may
  be viewed as a left $R^\nu$-module with respect to the action
  $$
  x *m:=\nu_{ij}xm
  $$
  for $x \in R_i$ and $m \in M_j$. To avoid confusion, we will write $M^\nu$ for $M$ viewed as an $R^\nu$-module and write $m^\nu$ to denote the element $m$ in $M$ viewed as an element in $M^\nu$.
  We now label the homomorphism $(\a,\alphaol)$ of graded rings as 
  $(\b,\betaol) :(R,\G) \to (R^\nu,\G)$.  As remarked at the end of section \ref{sect.gr.rings},
  there is an equivalence $\b_*:\Gr(R^\nu,\G) \to \Gr(R,\G)$ given by the following rule:  
   if $M$ is an $R^\nu$-module, then $\b_*M$ is an
  $R$-module with $x \in R_i$ acting on $\b_*m \in (\b_*M)_j=M_{\betaol j}$ by $x.\b_*m=\b_*(\b(x)m)$.

 If we now identify the domain of $\b_*$ with  $ \Gr(R^\nu,\G)$ of $\nu$ with $\Gr(R,\G)$ then $\b_*$ is 
 the auto-equivalence $\a_*$. To see this suppose that $M \in \Gr(R,\G)$ and consider 
 the action of $x \in R_i$ on an element $\b_*(m^\nu)$ in $\b_*(M^\nu)_j = (M^\nu)_{\betaol j} = 
 M_{\betaol j}=M_{\alphaol j}$. We have 
 $$
 x. \b_*(m^\nu) = \b_*(\b(x)*m^\nu)= \b_*(\nu_{\betaol i,\betaol j} (\b(x)m )^\nu) = \l_{ij} \b_*((\a(x)m)^\nu).
 $$
 Stripping away the superfluous notation this reads $x.m=\l_{ij}\a(x)m$ which is, indeed, the action of 
 $x$ on $\a_*M$.

 \subsubsection{Warning}
 \label{sect.warn1}
 Some care must be taken when identifying $\a_*M$ with its underlying set $M$, even when $\a$ is an automorphism. Suppose $c \in R_j$ 
 belongs to the center of $R$. Then the multiplication map $\rho_{c,M}:M \to M(j)$, $\rho_{c,M}(m):=mc$
 is a homomorphism of graded $R$-modules. By definition, $\a_*(\rho_{c,M})$ is the homomorphism 
 $\a_*M \to \a_*(M(j))$ given by 
 $$
 (\a_*(\rho_{c,M}))(\a_*m)= \a_*(\rho_{c,M}(m)) = \a_*(mc) = \a_*(m).\alpha^{-1}c  
 $$
  so
 $$
 \a_*(\rho_{c,M} )= \rho_{\a^{-1}c,\a_*M}.
 $$
 Thus $\a_*(\rho_{c,M})$ is multiplication by $\a^{-1} c$ on $\a_*M$ but when/if 
 $\a_*M$ is identified with its underlying set $M$, $\a_*(\rho_{c,M})$ is multiplication by $c$.  
 When $(\a,\l)$ is an almost automorphism $ \a_*(\rho_{c,M} )$ acts on $(\a_*M)_i$ as 
 multiplication by $\l_{ij}  \a^{-1} c$. 

\begin{lem}
\label{lem.theta}
Let $(\a,\l)$ be an almost-automorphism of $(R,\G)$ and  $h \in \G$. 
The map $\theta:R(h) \to \alpha_* (R(\alphaol h))$ defined by 
$$
\theta(x):= \l^{-1}_{h,i} \alpha_* (\a x) \qquad \hbox{if } \; x \in R(h)_i
$$ 
is an isomorphism of graded right $R$-modules. In particular, $\a$ viewed as a map $R \to \a_* R$ is an 
isomorphism of right $R$-modules.  
\end{lem}
\begin{pf}\footnote{This proof does not use the fact that $\l_{00}=1$.}
Let $x \in R( h)_i$. Then $x \in R_{ h+ i}$, so $\a x \in R_{\alphaol h+\alphaol i} = R(\alphaol h)_{\alphaol i}$.
Hence $\alpha_* (\a x) \in \alpha_* \bigl( R(\alphaol h)\bigr)_i$.  Thus $\theta$ preserves degree. 

Because $\a$ is bijective $\theta$ is too.

To see that $\theta$ is a right $R$-module homomorphism, suppose that $x\in R(h)_i$ and 
$y \in  R_j$. Then
\begin{align*}
\theta(x.y)  & =  \a_*(\a(xy))  \l_{h,i+j}^{-1}
\\
& = \alpha_* ( \a(x)\a(y) )  \l_{i+h, j}   \l_{h,i+j}^{-1} 
\\
&  = \a_*(\a (x)\a(y))   \l_{hi}^{-1} \l_{ij} 
\\
&=    \a_*(\a x).y    \l_{hi}^{-1}
\\
& =  \theta(x).y
\end{align*}
so $\theta$ is an $R$-module homomorphism. 
\end{pf}

 \begin{prop}
\label{prop.psi}
Let $(\a,\l)$ be an almost-automorphism of $(R,\G)$ and $h \in \G$. There is an isomorphism of functors
$$
 (h) \circ \a_* \cong \a_*  \circ  (\alphaol h).
 $$
 If $\a$ is an automorphism, then $ (h) \circ \a_*  =  \a_*  \circ  (\alphaol h)$.
\end{prop}
\begin{pf}
First we show there is an isomorphism of graded $R$-modules 
$$
 (\a_*M)(h) \cong \a_* (M(\alphaol h)) 
 $$
 for every $M \in \Gr(R,\G)$. 
 
Let $i \in \G$. 
The degree $i$ components of $\a_* (M(\alphaol h))$ and $(\a_*M)(h)$ 
are equal to $M_{\alphaol h + \alphaol i}$.
Let $m \in M_{\alphaol h+\alphaol i}$. For the purposes of this proof, we will write
 $\wtm$ for $m$ viewed as an  element of   $(\a_*M)(h)$ 
and $\whm$  for $m$ viewed as an element of  $\a_* (M(\alphaol h))$.
 
The map
$$
\psi:   (\a_*M)(h) \longrightarrow \a_* (M(\alphaol h)) , \quad \psi(\wtm) := \l_{h,i}^{-1} \whm
\quad \hbox{for } \; m \in (\a_*M)(h)_i
$$
preserves degree and is bijective because the $\l_{hi}$s are units. 
Furthermore, if $y \in  R_j$, then
\begin{align*}
\psi(\wtm.y)  & =  \psi\Bigl(\widetilde{m\a(y)}\Bigr) \l_{i+h,j}  
\\
& =  \widehat{m\a(y)}   \l_{h,i+j}^{-1}  \l_{i+h,j}  
\\
& =  \widehat{m\a(y)}   \l_{ij}  \l_{h,i}^{-1}  
\\
&=    \whm .y    \l_{hi}^{-1}
\\
& =  \psi(\wtm).y
\end{align*}
so $\psi$ is an $R$-module homomorphism, and hence an isomorphism of graded $R$-modules.

Let $f:M \to N$ be a homomorphism between graded $R$-modules. We write $\wtf$ and $\whf$
respectively for the homomorphisms obtained by applying the functors  
$(h) \circ \a_*$ and $\a_*  \circ  (\alphaol h)$ to $f$. Of course, $\wtf$ and $\whf$ are just the map $f$
on the underlying set $M$. 
It is easy to see that 
\begin{equation}
\label{eq.alpha*}
\UseComputerModernTips
  \xymatrix{
(\a_* M)(h) \ar[rr]^{\wtf}  \ar[d]_{\psi_M} && (\a_* N)(h)  \ar[d]^{\psi_N}
\\
 \a_* (M(\alphaol h)) \ar[rr]_{\whf} &&  \a_* (N(\alphaol h))
 }
\end{equation}
commutes. Hence the $\psi_M$s collectively give a natural isomorphism 
$ (h) \circ \a_* \to \a_*  \circ  (\alphaol h)$.

When $\a$ is an algebra automorphism we can take $\l_{hi} =1$ for all $h,i \in \G$, so the map 
$\psi$ is the identity map $\id_M$. Hence $ (h) \circ \a_* =  \a_*  \circ  (\alphaol h)$.
\end{pf}

 Lemma \ref{lem.theta} is a consequence of Proposition \ref{prop.psi} and the fact that $\a:R \to \a_* R$ 
 is an isomorphism of graded right $R$-modules.

 \subsection{Isomorphisms between auto-equivalences}
 Let $(R,\G)$ be a graded ring and $M$ a graded $R$-module.
 For each $y \in R_j$ and  each $i \in \G$, define 
 $$
 f_{i,y}:R(i)\to R(i+j), \qquad f_{i,y}(x):=yx.
 $$
 Then $f_{i,y}$ is a homomorphism of graded right $R$-modules.
  
 \begin{prop}
 \label{prop.wisom}
 Let $F$ and $G$ be auto-equivalences of  $\Gr(R,\G)$ such that 
  $$F(R(i)) \cong G(R(i))$$ 
 for  all $i \in \G$. 
 If the map $R_0 \to \End_{R_0}(R_j)$ sending $a
 \in R_0$ to the map $b \mapsto ab$ is an isomorphism for all $j \in \G$, then $$F \cong G.
 $$
  \end{prop}
 \begin{pf}
 If the proposition is true when $G=\id_{\Gr(R,\G)}$, then it holds for all $G$ because the general result 
 is obtained by applying the special case to $F'G$ where $F'$ is a quasi-inverse to $F$.
 We will therefore assume that $F(R(i)) \cong R(i)$ for all $i \in \G$ and show 
 that $F \cong \id_{\Gr(R,\G)}$. 
 
By hypothesis, there are isomorphisms $\phi_i:R(i) \to F(R(i))$ for all $i \in \G$.  

Fix $i \in \G$. Let $y \in R_j$. Then 
$$
\phi_{i+j}^{-1} \circ F(f_{i,y}) \circ \phi_i: R(i) \to R(i+j)
$$
is a homomorphism of graded right $R$-modules, so is left multiplication by a unique element $\theta_i(y) \in 
R_j$. That is,
$$
\phi_{i+j}^{-1} \circ F(f_{i,y}) \circ \phi_i = f_{i,\theta_i(y)}.
$$
We therefore have a map $\theta_i:R(j) \to R(j)$ for all  $j \in \G$.
If $a \in R_0$, then $f_{i,ay}= f_{i+j,a} \circ f_{i,y}$, so $f_{i,\theta_i(ay)} = f_{i,a\theta_i(y)}$. 
Hence $\theta_i$ is a right $R_0$-module homomorphism. 

If $w \in R_j$, then $\phi_{i+j} \circ f_{i,w} \circ \phi^{-1}_i:F(R(i)) \to F(R(i+j))$, but $F:\hom(R(i),R(i+j))
\to \hom(F(R(i)),F(R(i+j)))$ is bijective because $F$ is an auto-equivalence so 
$$\phi_{i+j} \circ f_{i,w} \circ \phi^{-1}_i = F(f_{i,y})$$ for some $y \in R_j$. Hence $f_{i,w}=f_{i,\theta_i(y)}$.
This proves that $\theta_i$ is bijective and hence an isomorphism of right $R_0$-modules.

In particular, $\theta_0:R(j) \to R(j)$ is an isomorphism of right $R_0$-modules for every $j$, so there
are units $u_j \in R_0$ such that $\theta_0(y)=u_jy$ for all $y \in R_j$.

Let $z \in R_k$. Then $f_{i+j,zy}=f_{k,z} \circ f_{i,y}$ so
\begin{equation}
\label{eq.wisom}
\theta_i(zy)=\theta_{i+j}(z) \theta_i(y).
\end{equation}
In particular, taking $i=0$ we see that
$$
u_{j+k} zy = \theta_j( z) u_j y
$$
for all $y \in R_j$ and all  $z \in R_k$.

For each $j \in \G$ define $\tau_j:R(j) \to F(R(j))$ to be $\tau_j: =u_j \phi_j$.
Let $z \in R_k$. 
Then the diagram
$$
\UseComputerModernTips
  \xymatrix{
R(j) \ar[rr]^{\tau_j}     \ar[d]_{f_{j,z}} && F(R(j)) \ar[d]^{F(f_{j,z})}
\\
R(j+k)   \ar[rr]_{\tau_{j+k}} && F(R(j+k))
 }
$$
commutes because if $y \in R(j)$, then
\begin{align*}
\bigl(F(f_{j,z}) \circ \tau_j \bigr)(y) & = \bigl( F(f_{j,z}) \circ u_j \phi_j \bigr)(y)
\\
 & = \bigl( F(f_{j,z}) \circ \phi_j \bigr) (u_j y)
\\
& = \bigl( \phi_{j+k} \circ f_{j,\theta_j(z)}\bigr)(u_j y)
\\
&=  \phi_{j+k} \bigl( \theta_j(z) u_j y\bigr)
\\
&=  \phi_{j+k} \bigl(u_{j+k} zy \bigr)
\\
& = \tau_{j+k}(zy)
\\
& = \bigl(\tau_{j+k} \circ f_{j,z}\bigr)(y).
\end{align*}
Thus, the $\tau_j$s taken together provide a natural isomorphism from the restriction of $F$
to the full subcategory of $\Gr(R,\G)$ consisting of the $R(j)$s to the identity functor on that
subcategory. 

Because the $R(j)$s generate $\Gr(R,\G)$ the isomorphism $\tau$ extends to an isomorphism $\tau:F \to
\id_{\Gr(R,\G)}$. (One defines $\tau_M:M \to FM$ for a general $M$ by writing 
$M$ as the cokernel of a map $P \to Q$ where $P$ and $Q$ are direct sums of various $R(j)s$.)
 \end{pf}

 The hypotheses of Proposition \ref{prop.wisom} hold   if $R_0$ is 
 an integrally closed commutative noetherian 
 domain and each $R_i$ is isomorphic to a non-zero ideal of $R_0$.
 
 \subsubsection{Proposition \ref{prop.wisom} applies to $A$ and $C$.}
Both $A_0$ and $C_\varnothing $ are isomorphic to a polynomial ring $k[z]$ and all the homogeneous components 
of $A$ and $C$ are rank one free $k[z]$-modules so the hypotheses of Proposition \ref{prop.wisom} are
satisfied by $A$ and $C$. 

Sierra uses the conclusion of Proposition \ref{prop.wisom}  for the Weyl algebra although her 
proof that the conclusion of 
 Proposition \ref{prop.wisom} holds for the Weyl algebra is rather different from our proof of 
 Proposition \ref{prop.wisom}.

\subsection{The Picard group $\Pic(R,\G)$}

There are several well-understood connections between $\Pic(R,\G)$ and other invariants of 
$(R,\G)$. For example, there is a group homomorphism $\G \to \Pic(R,\G)$ that sends $i \in \G$ to the twist
functor $M \mapsto M(i)$. The kernel of this map is the subgroup of $\G$ consisting of those $i$
for which $R_i$ contains a unit. The homogeneous units in $(A,\ZZ)$ and $(C,\ZZ_{\rm fin})$ belong to the 
identity component of the ring so the map  $\G \to \Pic(R,\G)$ is injective in those two cases.

The assignment $(\a,\alphaol,\l) \rightsquigarrow [\alpha_*]$ gives a map $\AAut \to \Pic(R,\G)$ and by Proposition \ref{prop.psi} the image of this map is contained in the normalizer of the image of $\G$ in 
$\Pic(R,\G)$.

 \subsubsection{Notation.} 
 If $\tau:F \to G$ is a natural transformation we write $\tau_M$ for the associated map $FM \to GM$. 

\begin{prop}
\label{prop.aut.end.id}
Let $\sA$ be an additive category.
There is group homomorphism 
$$
\Phi:\Pic(\sA) \to \Aut \bigl( \End (\id_{\sA} )\bigr)
$$
defined as follows:
If $F$ is an auto-equivalence of $\sA$, $G$ a left adjoint to $F$, and hence  a quasi-inverse to $F$, and 
$\eta:\id_{\sA} \to FG$ the unit,  
$$
\Phi([F]) :\End (\id_{\sA}) \to \End (\id_{\sA}), \qquad \Phi([F])(\tau)_M:=
 \eta^{-1}_M \circ F(\tau_{GM}) \circ \eta_M
 $$
\end{prop}
\begin{pf}
It is easy to check that $\Phi([F])$ is an automorphism of the ring $\End (\id_{\sA})$ and we omit the details. 
It is easy to check that $\Phi([F][F'])  = \Phi([F])\Phi([F'])$ {\it provided} $\Phi$ is well-defined.  To check
$\Phi$ is well-defined it suffices to show that $\Phi([F])$ is the identity if $F \cong \id_{\sA}$.

Let $\theta:F \to \id_\sA$ be an isomorphism. Let $\tau \in \End(\id_{\sA})$. 
The large rectangle in the diagram
$$
\UseComputerModernTips
  \xymatrix{
  M \ar[rr]^{\eta_M} \ar[dd]_{\tau_M}  && FGM \ar@{-->}[dd] \ar[rr]^{\theta_{GM}} && GM \ar[dd]^{\tau_{GM}}
  \\
  &&&&&&&
  \\
   M \ar[rr]_{\eta_M} && FGM \ar[rr]_{\theta_{GM}} & &  GM 
   }
   $$
commutes because $\tau$ is a natural transformation. If the dashed arrow is $F(\tau_{GM})$ then the right-hand square commutes because $\theta$ is a natural transformation $F \to \id_\sA$.
 If the dashed arrow is $FG(\tau_{M})$ then the left-hand square commutes because $\eta$ is a natural transformation $\id_\sA \to FG$. Because the horizontal arrows in the diagram are isomorphisms, it follows that $$F(\tau_{GM}) = FG(\tau_{M}).
 $$
  In particular,
 $$ 
 \Phi([F])(\tau)_M =  \eta^{-1}_M \circ F(\tau_{GM}) \circ \eta_M =  \eta^{-1}_M \circ FG(\tau_{M}) \circ \eta_M =\tau_M.
 $$
 This holds for all $M$ so $\Phi([F])(\tau) =\tau$. But $\tau$ was arbitrary, so $\Phi([F])$ is the identity. Hence
 $\Phi$ is a well-defined homomorphism. 
\end{pf}

One should check that the definition of $\Phi([F])$ does not depend on the choice of $G$. 

\smallskip

Suppose $R_0$ is central in $R$.
For each $a \in R_0$ and $M \in \Gr(R,\G)$, let $\mu_{a,M}:M \to M$ be  $\mu_{a,M}(m):=
ma$ for $m \in M$. Let $\mu_a:\id_{\Gr(R,\G)} \to id_{\Gr(R,\G)} $ be the natural transformation
$(\mu_a)_M:=\mu_{a,M}$.
Then the map 
\begin{equation}
\label{eq.R0.id}
\mu:R_0 \to \End\bigl(\id_{\Gr(R,\G)}\bigr), \qquad a \mapsto \mu_a. 
\end{equation}
 is a ring isomorphism.

\begin{prop}
\label{prop.aut.R0}
Let $(R,\G)$ be a graded ring and suppose that $R_0$ belongs to the center of $R$. 
Let $\Phi$ be the group homomorphism in Proposition \ref{prop.aut.end.id}, 
and let $\mu$ be the isomorphism in  (\ref{eq.R0.id}). Then the map
$$
\a: \Pic(R,\G) \to \Aut(R_0), \qquad [F] \mapsto \a_F := \mu^{-1} \circ \Phi([F]) \circ \mu
$$
is a group homomorphism and, if $M \in \Gr(R,\G)$, then $Ma=0$ if and only if $(FM).\a_F(a)=0$.
\end{prop}
\begin{pf}
Let $F$ be an auto-equivalence of $\Gr(R,\G)$, $G$ a left adjoint to $F$, and $\eta:\id_{\Gr(R,\G)} \to FG$ 
the unit.
 By definition, $\mu \circ \a_F=\Phi([F]) \circ \mu$, so  
 \begin{align*}
 (FM).\a_F(a)  & = {\rm Image}\bigl( \mu(\a_F(a))_{FM}   \bigr)
 \\
 &=  {\rm Image}\bigl( \Phi([F])(\mu_a)_{FM}   \bigr)
 \\
 & = {\rm Image}\bigl(  \eta^{-1}_{FM} \circ F(\mu_{a,GFM} ) \circ \eta_{FM}  \bigr).
 \end{align*}
 Hence $(FM).\a_F(a)=0$ if and only if $\mu_{a,GFM}=0$, if and only if $(GFM).a=0$, if and only if $Ma=0$.
\end{pf}

More succinctly,  $\a_F(a)$ is the unique $b \in R_0$   such that $\Phi([F])(\mu_a) = \mu_b$.

The next section examines a situation where there  is a homomorphism $R_0 \to \End(\id_{\Gr(R,\G)})$
without the hypothesis that $R_0$ is central in $R$. The result applies to $(A,\ZZ)$ and  is used
implicitly in the proof of \cite[Thm. 5.5]{Sue1}.

\subsection{$R_0$-$R$-bimodules}
\label{sect.R0.action}

Sierra exploits to advantage  the fact that every graded right 
$A$-module can be made into an $A_0$-$A$-bimodule. 
Proposition \ref{prop.R0}, the hypotheses of which are satisfied by $A$ and $C$,
gives a criterion on a graded ring $(R,\G)$ that implies every graded right 
$R$-module can be made into an $R_0$-$R$-bimodule. First we need a lemma.

\begin{lem}
\label{lem.R0}
Suppose that 
\begin{enumerate}
  \item 
  $R_0$ is commutative;
  \item 
   $R_i$ is a torsion-free left $R_0$-module for all $i \in \G$;
  \item 
  $R_0x=xR_0$ for all $x \in R_i$ and all $i \in \G$;
  \item
  $R_0x \cap R_0 y \ne 0$ for all $x,y \in R_i-\{0\}$ and all $i \in \G$.
\end{enumerate}
Then there is a homomorphism $\G \to \Aut(R_0)$, $i \mapsto \theta_i$, such that $xa=\theta_i(a)x$
for all $a \in R_0$, all $x \in R_i$, and all $i \in \G$.
\end{lem}
\begin{pf}
Fix $i \in \G$. 
Let $a \in R_0$.

\underline{Claim:} 
There is a unique $a' \in R$ such that $xa=a'x$ for all $x \in R_i$.

\underline{Proof:}
Let $x,y \in R_i-\{0\}$. By hypothesis (2), there are elements $a',a'' \in R_0$ such that 
$xa=a'x$ and $ya=a'' y$. 
By hypothesis (4), there are $b,c \in R_0$ such that $bx=cy \ne 0$.
Therefore
$$
a'bx=ba'x=bxa=cya=ca''y=a''cy=a''bx,
$$
and $(a'-a'')bx=0$. But $bx \ne 0$ so hypothesis (2) implies that $a'=a''$. $\lozenge$

It follows that there is a well-defined map $\theta_i:R_0 \to R_0$ such that $xa=\theta_i(a)x$
for all $a \in R_0$ and all $x \in R_i$. Let $b \in R_0$. It is clear that $\theta_i(a+b)=\theta_i(a)+\theta_i(b)$.
Also, $\theta_i(ab)x=xab=\theta_i(a)xb=\theta_i(a) \theta_i(b)x$ so, by hypothesis (2), $\theta_i(ab)=
 \theta_i(a) \theta_i(b)$. Hence $\theta_i$ is an endomorphism of $R_0$. By hypothesis (3),
 $\theta_i$ is surjective. By hypothesis (2), $\theta_i$ is injective. Hence $\theta_i \in \Aut(R_0)$.
 
Let $j \in \G$. Let $a \in R_0$, $x \in R_i$, and $y \in R_j$. Then 
$$
\theta_{i+j}(a)xy=xya=x\theta_j(a) y =\theta_i\theta_j(a) xy
$$
so, by hypothesis (2), $\theta_{i+j}  =\theta_i\theta_j$. Hence $i \mapsto \theta_i$ is a group
homomorphism, as claimed. 
\end{pf}

\begin{prop}
\label{prop.R0}
Suppose that $(R,\G)$ satisfies the hypotheses in Lemma \ref{lem.R0} and let 
$\theta_i \in \Aut(R_0)$, $i \in \G$, be defined as in loc. cit.. 
Then 
\begin{enumerate}
  \item 
  every $M \in \Gr(R,\G)$ is an $R_0$-$R$-bimodule with respect to the action
$$
a.m:= m\theta_{-i}(a) \qquad \hbox{for $m \in M_i$ and $a \in R_0$};
$$
  \item 
  the map 
  $$\mu:R_0 \to \End(\id_{\Gr(R,\G)}), \quad \mu(a)_M:M \to M, \quad \mu(a)_M(m):= a.m,
  $$
   is an isomorphism of rings.
\end{enumerate}
\end{prop}
\begin{pf}
(1)
If $a,b \in R_0$, then 
$$
(ba).m=(ab).m= m\theta_{-i}(ab)=m\theta_{-i}(a)\theta_{-j}(b) = (a.m)\theta_{-i}(b) = b.(a.m)
$$
so $M$ is a left $R_0$-module.
If $y \in R_j$, then
$$
a.(my)=my\theta_{-i-j}(a)=my\theta_{-j}\theta_{-i}(a)=m\theta_{-i}(a)y=(a.m)y.
$$
Hence $M$ is an $R_0$-$R$-bimodule.

(2)
Because left multiplication by $a \in R_0$ is an $R$-module endomorphism of $M$, $\mu(a)$ is  
a natural transformation. It is easy to check that $\mu$ is a homomorphism. It is obviously injective. 

To see that $\mu$ is surjective, let $\tau:  \id_{\Gr(R,\G)} \to \id_{\Gr(R,\G)}$ be a natural transformation. 
Define $a :=\tau_R(1)$. 
Let $r \in R_i$, and let 
 $\l_r:R(-i) \to R$ be the map $\l_r(x)=rx$. The diagram
$$
\UseComputerModernTips
  \xymatrix{
R(-i) \ar[r]^{\l_r} \ar[d]_{\tau_{R(-i)}} & R     \ar[d]^{\tau_{R}}
\\
R(-i) \ar[r]_{\l_r}  & R  
}
$$
commutes so 
$$
r\tau_{R(-i)}(1) = \tau_R\l_r(1) = \tau_R(r)=\tau_R(1.r) = ar = r \theta_{-i}(a).
$$
Hence $\tau_{R(-i)}(1) =  \theta_{-i}(a)$. Now consider an arbitrary graded right $R$-module $M$ 
and an element $m \in M_i$.
Let $\l_m:R(-i) \to M$ be the map $f(x)=mx$. The diagram 
$$
\UseComputerModernTips
  \xymatrix{
R(-i) \ar[r]^{\l_m} \ar[d]_{\tau_{R(-i)}} & M     \ar[d]^{\tau_{M}}
\\
R(-i) \ar[r]_{\l_m}  & M  
}
$$
commutes so 
$$
a.m = m \theta_{-i}(a) = \l_m \bigl(\tau_{R(-i)}(1) \bigr) = \tau_M \l_m(1)=\tau_M(m) 
$$
for all $m \in M$. It follows that $\tau_M = \mu(a)_M$ and that $\tau=\mu(a)$. Hence $\mu$ is surjective. 
\end{pf}

\subsection{}
There is one further way in which an auto-equivalence of $\Gr(R,\G)$ can induce an automorphism
of $R_0$. The map 
$$
\l_j:R_0 \to \hom(R(j),R(j)), \qquad    \l_j(a)(x):=ax.
$$
 is an isomorphism of rings. The following result is therefore clear.

\begin{lem}
Let $F$ be an autoequivalence of $\Gr(R,\G)$.
Let $f:\hom(R,R) \to \hom(FR,FR)$ be the isomorphism 
$g \mapsto Fg$. If $FR \cong R(j)$, then   $\l_j^{-1} \circ f \circ \l_0$ is an automorphism of $R_0$. 
\end{lem}

\section{$C$ as an ungraded ring}
\label{sect.C.ring}

In this section, as in others, we assume that $k$ is of characteristic zero.

The results in this section are not required for the proof of the main result in the paper
but  $C$ is an interesting example of a class of commutative rings not commonly encountered 
so we establish its basic properties here.

\subsection{}
If $I \subset \ZZ-\{0\}$, we write $R_I$ for the subring of $C$ generated by $\{x_0\} \cup \{x_n \; | \;  n \in I\}$.

\begin{prop}
\label{prop.DD}
The ring $C$ is an ascending union of Dedekind domains, and is flat over each of those 
Dedekind domains.
\end{prop}
\begin{pf}
It is clear that $C$ is the ascending union of the subrings $R_I$ where the union is taken
over any ascending and exhaustive chain of finite subsets $I\subset \ZZ-\{0\}$. Such a subring is  isomorphic to the ring
$$
S_I:= {{k[t,X_n\; | \; n \in I]}\over{\fa}}
$$
where $\fa$ is the ideal generated by the elements
$$
g_n:=X_n^2-t^2+n, \quad n \in I.
$$
(The element $t$ corresponds to $x_0$.)

Let $I$ be a finite subset of $\ZZ-\{0\}$. Let $k(t)$ be the rational function field over $k$ and 
let $F$ be a splitting field for the polynomial 
$$
f(X)=\prod_{n \in I} (X^2-t^2+n).
$$
Constructing $F$ as a tower of quadratic extensions, it is easy to see that the integral closure of $k(t)$ in
$F$ is isomorphic to $S_I$, and hence to $R_I$. Therefore $R_I$ is a Dedekind domain.

If $I \subset J \subset \ZZ-\{0\}$ are finite subsets, then $R_I$ is contained in $R_J$ and $R_J$ is a finitely generated torsion-free, and hence projective, $R_I$-module. Hence, for every $I$, $C$ is a directed 
union of finitely generated projective $R_I$-modules, and is therefore a flat $R_I$-module.
\end{pf}
 
\subsubsection{}

There are other ways to prove  Proposition \ref{prop.DD}. For example, one can prove directly, using the Jacobian criterion, that the rings $S_I$ in its proof are regular of Krull dimension one. 

\subsubsection{}
Suppose $k$ is algebraically closed and  fix elements $\sqrt{-n}$ in $k$.
Let $k[[z]]$ be the ring of formal power series.
There is a homomorphism $\varphi:C \to k[[z]]$ given by
$$
\varphi(x_0):=z, \qquad \varphi(x_n):=   \sqrt{-n}\biggl( 1- {{z^2}\over{n}}\biggr)^{1/2}, \quad n \ne 0,
$$
where $(1-z^2/n)^{1/2}$ denotes the Taylor series expansion for $\sqrt{1-z^2/n}$ centered at $z=0$.
The restriction of $\varphi$ to the Dedekind domains $S_I$ appearing  in the proof of Proposition \ref{prop.DD} is injective, so $\varphi$ is injective on $C$.

\begin{prop}
The ring $C$ has the following properties:
\begin{enumerate}
  \item 
 It is an  integrally closed  non-noetherian domain.
\item{}
Its transcendence degree is one.
\item
Suppose $k = \CC$. 
If $\fm$ is a maximal ideal in $C$, then $\dim_k(\fm/\fm^2)=1$.  
\item
Every finitely generated ideal in $C$ is projective and generated by $\le 2$ elements.
  \item 
  Let $d$ be a positive integer. The ring homomorphism $\g:C \to C$ defined by $\c(x_n):=x_{nd}/\sqrt{d}$
is an isomorphism from $C$ onto its subalgebra $k[x_{nd} \; | \; n \in \ZZ]$.
\end{enumerate}
\end{prop}
\begin{pf}
(1)
Of course, $C$ is a domain because it is an ascending union of domains. It is integrally closed because it
is a directed union of integrally closed rings.

To show $C$ is not noetherian it suffices to show that $C \otimes_k \kol$ is not noetherian so we may,
and will,  assume $k$ is algebraically closed. 

For each integer $N$, let $\fa_N$ be the ideal generated by the elements $x_{d}+\sqrt{-d}$ for $d \le N$.
Then $\fa_N \subset \fa_{N+1}$ but   $\fa_N \ne\fa_{N+1}$.

(2)
Since the field $\Fract C$ is the union of finite extensions of $k(x_0)$, it is clear that 
$x_0$ is a transcendence basis for $\Fract C$.

(3)
Because $C$ has countable dimension whereas the rational function field $\CC(t)$ has uncountable dimension,  
$\fm$ is generated by $\{x_n-z_n \; | \; n \in \ZZ\}$ for suitable elements $z_n \in \CC$. 
We write $z$ for the point $(z_n)_{n \in \ZZ} \in \CC^\ZZ$ and think of it as a closed point of $\Spec C$.

The same argument as for the polynomial ring in a finite number of variables
shows that  $\fm=\fm^2+\sum_n k(x_n-z_n)$.

Fix an integer $r$.
Because $z \in \Spec C$, $z_n^2+n=z_r^2+r$ for all $n$, whence
$$
z_n^2 -z_r^2=r-n=x_n^2-x_r^2.
$$
Therefore $\fm^2+k(x_r-z_r)$ contains
$$
\hbox{${{1}\over{2}}$}\biggl((x_r-z_r)^2 - (x_n-z_n)^2\biggr) + z_r(x_r-z_r)= z_n(x_n-z_n).
$$
If $z_n \ne 0$, then $x_n-z_n \in \fm$. 
In particular, if all $z_n$ are non-zero,  then $\fm^2+k(x_r-z_r)=\fm$. 

On the other hand, suppose $z_r=0$. Then $z_n \ne 0$ for all $n \ne r$
so the argument just given shows that $x_n-z_n \in \fm^2+kx_r$. Of course, $x_r \in \fm^2+kx_r$ too, and therefore $\fm^2+kx_r=\fm$.

(4)
A finitely generated ideal in $C$ is generated by an ideal in $R_I$ for some finite subset $I \subset \ZZ-\{0\}$.   However, every ideal in $R_I$ is projective. If $\fa$ is an ideal in $R=R_I$, then the kernel of 
the multiplication map $C \otimes_R \fa \to C\fa$ is isomorphic to $\Tor^R_1(C,R/\fa)$ which is zero because 
$C$ is flat over $R$. Thus $C\fa \cong C \otimes_R \fa$ and hence $\fa$ is a projective $C$-module.

Every ideal in a Dedekind domain can be generated by $\le 2$ elements.

(5)
This is straightforward.
\end{pf}

 \subsection{}
 \label{sect.RS}
 The results in this section will not be used elsewhere in this paper.
 
Let $k=\CC$, and give $\CC$ its usual topology. 

Let's write $\CC^\ZZ$ for the linear dual of the vector space $\oplus_{n \in \ZZ} \CC x_n$ and
 view the $x_n$s as coordinate functions on $\CC^\ZZ$ .  
 The $\CC$-valued points of $\Spec C$ are  the points in
  $$
  X:=\{ z=(z_n)_{n \in \ZZ} \; | \; \hbox{$z_n^2=z_0^2-n$ for all $n \in \ZZ$}\} \, \subset \, \CC^\ZZ.
  $$
 We now consider the question of whether $X$ can be given the structure of a Riemann surface.  
  In order to preserve the usual
 connection between complex algebraic curves and Riemann surfaces, we are particularly interested in 
 whether $X$ can be made into a Riemann surface in such a way that the coordinate functions 
 $x_n$ are holomorphic.  We will show this cannot be done when $\CC^\ZZ$ is given the product
 or the box topologies.
 On the other hand, if $\CC^\ZZ$ is identified with $\ell^\infty(\ZZ)$ in a suitable way, then $X$ has uncountably
 many connected components, all homeomorphic to one another, and each component can be 
 given the structure of a Riemann surface in such a way that each $x_n$ is a holomorphic function. 
 
 I thank Robin Graham for telling me the following result and allowing me to include it here.
 
\begin{prop}
\label{prop.box}
If $\CC^\ZZ$ is given the box topology, then $X$ is discrete.
\end{prop}
\begin{pf}
For $t \in \CC$ and $r \in \RR_{>0}$, let $D(t,r)$ denote the open disk of radius $r$ centered at $t$.
Fix a point $z=(z_n)$ in $X$. If $z_n=0$, let $r_n:={{1}\over{2|n|+1}}$.
If $z_n \ne 0$, let
$$
r_n:=\min\bigg\{|z_n|, {{1}\over{2|n|+1}} \bigg\}.
$$
Since $|z_n| \to \infty$ as $|n| \to \infty$, $r_n=|z_n|$ for only finitely many $n$. 
Since 
$$
U := \prod_{n \in \ZZ}  D(z_n,r_n)
$$
is an open neighborhood of $z$ in the box topology, to show $X$ is discrete it suffices to show that
$U \cap X = \{z\}$.

Suppose  that $y=(y_n) \in X-\{z\}$. There are two cases.

(1)
Suppose $y_0^2=z_0^2$. Then $y_n =-z_n \ne 0$ for some $n$. 
For that $n$,  $|y_n-z_n|=2|z_n|$ so $y_n \notin D(z_n,r_n)$. Hence $y \notin U$.

(2)
Suppose $y_0^2 \ne z_0^2$. 
If $a,b \in \CC$ are fixed and $w \in \CC$, then 
$$
 \bigg\vert {{\sqrt{a+w} + \sqrt{b+w}  }\over{w}}  \, \bigg\vert  \longrightarrow 0
$$
regardless of which branches of the square root function are chosen, and regardless of whether 
the branches chosen for $\sqrt{a+w}$ and $ \sqrt{b+w}$ are the same or not.
Therefore, if $|n| \gg 0$, 
\begin{align*}
|y_n-z_n|  & = \bigg\vert \sqrt{y_0^2-n} - \sqrt{z_0^2-n}\, \bigg\vert 
\\
 & =  \bigg\vert   {{y_0^2-z_0^2}\over{\sqrt{y_0^2-n} + \sqrt{z_0^2-n}}} \, \bigg\vert 
\\
& > {{1}\over{|n|}}
\\
&  > {{1}\over{2|n|+1}}.
\end{align*}
Hence $y \notin U$, and $X$ is discrete, as claimed.
\end{pf}

I thank Lee Stout for telling me the following result and allowing me to include it here. 
  
 \begin{prop}
 \label{prop.not.RS}
 Let $\CC$ have its usual topology, $\CC^\ZZ$ the product topology, and $X$ the subspace topology.
 \begin{enumerate}
  \item 
  Every fiber of every $x_k$ is homeomorphic to a Cantor set.
  \item 
  $X$ cannot be made into a complex manifold in such a way that any of the coordinate functions is holomorphic.
\end{enumerate} 
 \end{prop}
 \begin{pf}
  Fix a point $y =(y_n) \in X$ and an integer $k$.
 Let $x_k:X \to \CC$ be the function taking the $k^{\th}$ coordinate. 
 
 (1)
  Let $F$ be the fiber of $x_k$ over the point $y_k$. 
 Then
 $$
 F=\prod_{n < k} \{-y_n,y_n\} \times \{y_k\} \times \prod_{n >k} \{-y_n,y_n\}
 $$
However, at most one $y_n$ is zero so $F$ is homeomorphic to a countable product 
of copies of the discrete  space $\{\pm\}$ endowed with the product topology. 
Therefore $F$ is a Cantor set.   
 (2)
 Let $U$ be any open neighborhood of $y$. By shrinking $U$ we can assume 
 there is a positive integer  $N>|k|$ and $\ve >0$ such that  
 $$
U=X \cap  \{(z_n)  \; | \; |z_n-y_n| <  \ve  \; \hbox{ if } \, -N \le n \le N\}.
$$
Then $U$ contains a point $z$ such that $z_k \ne y_k$. 
For  each integer $n$ in $[-N,N]$, choose the branch of the square root function such that
$y_n=\sqrt{x_k^2+k-n}$ and define $s_n(\xi)=\sqrt{\xi^2+k-n}$  for $\xi$ in a sufficiently small disk centered at
$y_k$.

Hence  $U$ contains the set
$$
F':=\{(x_n)  \; | \;   x_n=y_n \; \hbox{ if } \, -N \le n \le N\}.
$$
This is also a Cantor set, and $z_k-y_k$ vanishes on it. But every point of $F'$ is a limit point in $F'$.
An analytic function on an open set $U$ that vanishes on a subset having a limit point
is identically zero in a neighborhood of that limit point  (the coefficients in the Taylor series expansion
around the limit point are zero). Hence $z_k -y_k$ vanishes on $U$. 
But that is absurd because $U$ contains a point $(x_n)$ with $x_k \ne y_k$.
 \end{pf}

 The space of doubly infinite $\CC$-valued sequences $\xi=(\xi_n)$ such that
 $|\xi_n|$ is bounded is denoted by $\ell^\infty(\ZZ)$. It is a Banach space with respect to the norm
 $$
 \parallel \xi \parallel = \sup_n\{\xi_n\}.
 $$
   Let $f:X \to \ell^\infty(\ZZ)$ be the map
   $$
     f(z)= \biggl( {{z_n} \over {2n+1}} \biggr).
   $$
   Write $Y:=f(X)$. 
   If $y_n$, $n \in \ZZ$, are the obvious coordinate functions on $\ell^\infty(\ZZ)$, then $Y$ is the locus cut out
   by    the equations
   \begin{equation}
   \label{eq.l.infty}
    (1+2n)^2y_n^2 = y_0^2-n.
   \end{equation}
   The ring $\CC[y_n \; | \; n \in \ZZ]$ where the $y_n$s satisfy the relations (\ref{eq.l.infty}) is isomorphic
   to $C$. 
   
  We write $\{\pm 1\}^{(\ZZ)}$ for the subgroup of  $\{\pm 1\}^\ZZ$ consisting of the 
  functions $\ZZ \to \{\pm 1\}$ that take the value $-1$ only finitely often.
  It is a countable direct sum of copies $\{\pm 1\}$.
   
   The next result is due to Robin Graham. I am grateful for his allowing me to include it here. 
   
   \begin{prop}
   \label{prop.RS}
   Let $\ell^\infty(\ZZ)$ have its usual topology, and give $Y$ the subspace topology.
   Then 
   \begin{enumerate}
  \item 
  $Y$ has uncountably many connected components, 
  \item 
  all those components are homeomorphic to one another,
  \item 
  they are permuted transitively by the action of $\{\pm 1\}^\ZZ$, 
  \item
  each component is stable under the action of $\{\pm 1\}^{(\ZZ)}$, and
  \item{}
  each component can be given the structure of a Riemann surface in such a way that
  $C$ consists of holomorphic functions.
\end{enumerate}
   \end{prop}

\section{The $\ZZ_{\rm fin}$ graded ring $C$}
\label{sect.C}
\label{sect.C.gr}
     
In this section we establish the basic properties of $C$ {\it as a graded ring}. 
One of the main results is that every graded ideal
of $C$ is principal. Because $C$ is also a domain the standard results about modules over a PID 
carry over to the category of graded $C$-modules. In particular, every projective graded $C$-module
is a direct sum of twists of $C$. 

We end the section with the proof that  $\Gr A \equiv \Gr(C,\ZZ_{\rm fin})$. 
This is done by exhibiting a bigraded $A$-$C$-bimodule that is, as a $C$-module, a projective 
generator in $\Gr(C,\ZZ_{\rm fin})$and has $A$ as its graded endomorphism ring. 

Most questions about $C$ reduce to combinatorial questions about $\ZZ_{\rm fin}$.

\subsection{}
The homogeneous components of $C$ are 
$$
C_J=C_\varnothing  x_J.
$$

\begin{lemma}
Let $I,J,I',J' \in \ZZ_{\rm fin}$. Then
\begin{enumerate}
  \item 
  $x_Ix_J=x^2_{I \cap J} x_{I \oplus J}$;
  \item 
  $(Cx_I)_J= C_\varnothing  x^2_{I-J} x_J$;
  \item 
  the following conditions are equivalent:
  \begin{enumerate}
  \item 
  $x_Ix_J=x_{I'}x_{J'}$;
  \item 
  $I \cap J=I'\cap J'$ and $I \oplus J=I' \oplus J'$;
  \item 
   $I \cap J=I'\cap J'$ and $I \cup J=I' \cup J'$.
\end{enumerate}
\end{enumerate}
\end{lemma}
\begin{pf}
(1)
This follows from the identity $I\cup J = (I \cap J) \sqcup (I\oplus J)$.

(2)
We have  $(Cx_I)_J=C_{I \oplus J} x_I = C_\varnothing  x_{I \oplus J}x_I = C_\varnothing  x^2_{I-J} x_J$.

(3)
The equivalence of (b) and (c) follows from the identity
$$
I \cup J=I \oplus J \oplus (I \cap J).
$$
The equivalence of (a) and (b) follows from (1).
\end{pf}

 The identity component of 
the ring obtained by inverting all non-zero homogeneous elements of $C$ is the field of fractions of $C_\varnothing $,
the rational function field $k(x_0^2)$.

\subsection{Graded ideals in $C$}

Because each $x_I$ is a regular element of $C$,  
$$
Cx_I \cong C(I).
$$
Here $C(I)$ is the degree-shifted module: $C$ viewed as a graded module with 1 placed in degree $I(=-I)$. 

\begin{lem}
\label{lem.K0}
Let $I,J \in \ZZ_{\rm fin}$. Then
\begin{enumerate}
  \item 
  $Cx_I+Cx_J=Cx_{I \cap J}$ and 
  \item 
  $Cx_I \cap Cx_J=Cx_{I \cup J}$. 
\end{enumerate}
\end{lem}
\begin{pf}
(1)
Since $I=(I-J) \sqcup (I \cap J)$ and $J = (J-I) \sqcup (I \cap J)$, we have
$$
Cx_I+Cx_J =(Cx_{I-J} + Cx_{J-I})x_{I\cap J}.
$$
However, $I-J$ and $J-I$ are disjoint so $x_{I-J}^2$ and $x_{J-I}^2$ are relatively prime elements of 
$C_\varnothing $ whence $Cx_{I-J} + Cx_{J-I} = C$. Hence $Cx_I+Cx_J=Cx_{I \cap J}$.

(2)
Since 
\begin{align*}
(Cx_I \cap Cx_J)_K  & = C_\varnothing  x^2_{I-K}x_K \cap  C_\varnothing  x^2_{J-K}x_K 
\\
& =  C_\varnothing  x^2_{(I-K) \cup (J-K)}x_K 
\\
& = C_\varnothing  x^2_{(I \cup J)-K}x_K
\\
& = (Cx_{I \cup J})_K
\end{align*}
 for all $K \in \ZZ_{\rm fin}$,   $Cx_I \cap Cx_J=Cx_{I \cup J}$ as claimed.
 \end{pf}

\begin{prop}
\label{prop.gr.pid}
Every graded ideal of $C$ is generated by a single homogeneous element.
\end{prop}
\begin{pf}
Let $\fa$ be a non-zero graded ideal of $C$.
Let $d$ be a non-zero element of $C_\varnothing =k[z]$ of minimal $z$-degree with the property that 
$dx_I \in \fa$ for some $I \in \ZZ_{\rm fin}$. Let $J$ be of minimal cardinality such that $dx_J \in \fa$.

Let $fx_I$ be an arbitrary element of $\fa$ with $0 \ne f \in C_\varnothing $.  Let $h$ be the greatest common divisor of 
$d$ and $f$ in $C_\varnothing $. Both $x_I$ and $x_J$ divide $x_{I \cup J}$ 
so $dx_{I\cup J}$ and $f x_{I \cup J}$ belong to $\fa$. Therefore  $hx_{I \cup J}$ belongs to $\fa$.
 But $\deg h \le \deg d$ so the choice of $d$ implies that $\deg h=\deg d$. Therefore $d$ divides $f$.

Write $f=dg$ where $g \in k[z]$. Then 
$$
(dx_J,fx_I)=(dx_{J-I},fx_{I-J})x_{I \cap J} = (x_{J-I},gx_{I-J})dx_{I \cap J}.
$$
However, $x_i$ is a unit modulo $x_j$ if $i \ne j$, so $x_{I-J}$ is a unit modulo $x_{J-I}$. Hence
$$
(dx_J,fx_I)=  (x_{J-I},g)dx_{I \cap J}.
$$
But
$$
(x_j,g)=
\begin{cases}
	(x_j) & \text{if $x_j^2 | g$,}
	\\
	C & \text{otherwise,}
\end{cases}
$$
so 
$$
 (x_{J-I},g) = (x_K)  \qquad \hbox{where } \;  K=\{j \in J-I \; | \; x_j^2 \hbox{ divides } g\}.
 $$
 Therefore
 $$
(dx_J,fx_I)=  (dx_{I \cap J}x_K) = (dx_L)
$$
where $L=(I \cap J) \cup K \subset J$.  By the choice of $J$, the cardinality of $L$ can be no smaller than that
of $J$ so $L=J$ and $(dx_J,fx_I) = (dx_J)$. It follows that $\fa=(dx_J)$. 
\end{pf}

\begin{prop}
\label{prop.K0}
Let $I,J,I',J' \in \ZZ_{\rm fin}$. 
\begin{enumerate}
\item
There is an isomorphism of graded $C$-modules
\begin{equation}
\label{eq.K0}
Cx_I \oplus Cx_J \cong Cx_{I \cup J} \oplus Cx_{I \cap J}.
\end{equation}
  \item 
There is a surjective degree zero $C$-module homomorphism $$Cx_I \oplus Cx_J \to Cx_K$$ if and only if $I \cap J \subset K \subset I \cup J$.
    \item 
 There is an isomorphism of graded $C$-modules
  $$Cx_I \oplus Cx_J \cong Cx_{I'} \oplus Cx_{J'}$$ if and only if $I \cup J=I' \cup J'$ and
$I \cap J=I'\cap J'$. 
\end{enumerate}
\end{prop}
\begin{pf}
(1)
By Lemma \ref{lem.K0}, the exact sequence 
$$
0 \to Cx_I \cap Cx_J \to Cx_I \oplus Cx_J \to Cx_I + Cx_J \to 0
$$
can be rewritten as  
$$
0 \to Cx_{I \cup J}  \to Cx_I \oplus Cx_J \to Cx_{I \cap J} \to 0.
$$
The right-most term is projective, so the sequence splits giving the claimed isomorphism.

(2)
($\Leftarrow$)
This follows from (1) because the hypothesis on $K$ implies there is a set $L \in \ZZ_{\rm fin}$
 such that $K\cap L=I \cap J$ and $K \cup L=I\cup J$, namely $L=(I \cup J - K) \cup(I \cap J)$. 

($\Rightarrow)$
Suppose there is a surjective  degree zero $C$-module homomorphism  $Cx_I \oplus Cx_J \to Cx_K$.
Because $Cx_I \cong C(I)$, and so on, 
there is a surjective  degree zero $C$-module homomorphism   
$C(I \oplus K) \oplus C(J \oplus K) \to C$ and hence a surjective  degree zero $C$-module homomorphism $f:Cx_{I \oplus K} \oplus Cx_{J\oplus K} \to C$.
Since $f$ is completely determined by $f(x_{I \oplus K},0)$ and $f(0,x_{J \oplus K})$ which 
must belong to  $C_{I \oplus K}$ and $C_{J \oplus K}$ respectively, i.e., to $C_\varnothing  x_{I \oplus K}$ and  
$C_\varnothing  x_{J \oplus K}$, the image of 
$f$ is contained in $Cx_{I \oplus K} + Cx_{J\oplus K}$ which is equal to 
$Cx_{(I \oplus K) \cap (J \oplus K)}$ by Lemma \ref{lem.K0}. 
Hence $Cx_{(I \oplus K) \cap (J \oplus K)} =C$. Therefore $(I \oplus K) \cap (J \oplus K)=\varnothing$
and this implies that  $I \cap J \subset K \subset I \cup J$. 

(3)
($\Leftarrow$)
This follows from (\ref{eq.K0}).

($\Rightarrow)$
Suppose that $Cx_I \oplus Cx_J \cong Cx_{I'} \oplus Cx_{J'}$.
Because $Cx_{I'}$ and $Cx_{J'}$ are quotients of $Cx_I \oplus Cx_J$, (2) implies that 
$I \cap J \subset I' \subset I \cup J $ and $I \cap J \subset J' \subset I \cup J $, i.e., $I \cap J \subset I' \cap J'$. and $I' \cup J' \subset I \cup J$.
By symmetry, the reverse inclusions also hold.\footnote{Although the equalities $I \cap J=I'\cap J'$ and 
$I \cup J=I'\cup J'$ imply that $I \oplus J = I' \oplus J'$, this latter equality follows directly from  the isomorphism $Cx_I \oplus Cx_J \cong Cx_{I'} \oplus Cx_{J'}$ because
taking the second exterior power implies that $Cx_Ix_J \cong Cx_{I'}x_{J'}$ whence $\deg x_Ix_J=\deg x_{I'}x_{J'}$ and $I \oplus J = I' \oplus J'$.}
\end{pf}

\subsection{Torsion-free  and projective graded $C$-modules}

A graded module over a graded ring is said to be a 
\begin{enumerate}
  \item 
  {\sf free graded module} if it has a basis consisting of homogeneous elements;
  \item 
  {\sf projective graded module} if it is projective as an object in $\Gr(C,\ZZ_{\rm fin})$.
\end{enumerate}  

Let $M$ be a graded module over a commutative graded ring $R$. A homogeneous 
element $m \in M$ is 
{\sf torsion} if $xm=0$ for some homogeneous regular element $x \in R$. A graded module is 
{\sf torsion} if every homogeneous element in it is torsion and  {\sf torsion-free} if
its only torsion element is 0. The submodule of $M$ 
generated by the torsion elements is a torsion module and is called the {\sf torsion submodule}
of $M$. We will denote it by $\tau M$ for now. The quotient $M/\tau M$ is torsion-free.

Presumably the following result is already in the literature.
     
     \begin{prop}
     \label{prop.pid.mods}
 Let $(R,\G)$ be a graded ring. Suppose $R$ is commutative, that all homogeneous elements of $R$
 are regular, and that every graded ideal  of $R$ is principal. Then 
\begin{enumerate}
  \item 
  every graded submodule of a finitely generated free graded module is a free graded module;\footnote{The finitely generated hypothesis can be removed.}
  \item 
   every finitely generated graded $R$-module is a direct sum of a graded torsion module and a 
 free  graded module.
 \end{enumerate}
 \end{prop}
 \begin{pf}
 (1)
 Let $f_1,\ldots,f_n$ be a homogeneous basis for a graded module $F$. We argue by induction on $n$
 to show that every graded submodule of $F$ is a free graded module.
 The result is true by hypothesis  if $n=1$ so suppose $n \ge 2$.
 
 Let $E$ be a graded submodule of $F$. 
 Let $\a:F \to Rf_1$ be the projection with kernel $F':=Rf_2 \oplus \cdots \oplus Rf_n$. If $E \subset F'$
 then $E$ has a homogeneous basis by the induction hypothesis, so we may suppose that $E \not\subset F'$.
 Then $\a(E)$ is a non-zero graded submodule of $Rf_n$, so is equal to $Raf_n$ for some homogeneous
 $a \in R$. But $Raf_n$ is isomorphic to a twist of $R$ so the map $\a|_E:E \to \a(E)$ splits and $E \cong \a(E)
 \oplus (E \cap F')$. By the induction hypothesis, $E \cap F'$ has a homogeneous basis. Hence $E$ has
 a homogeneous basis. 
 
 (2)
 Let $M$ be a finitely generated graded $R$-module. Since $M\tau M$ is torsion-free it suffices to show
 that a torsion-free finitely generated graded $R$-module is a free graded module.  So, we assume $M$
 is torsion-free.
 
 Let $\{m_1,\ldots,m_n\}$ be a homogeneous set of generators of $M$ and assume they have been 
 ordered so that $\{m_1,\ldots,m_s\}$ is a maximal subset of linearly independent elements. Write 
 $L=Rm_1+ \cdots + Rm_s$. If $s=n$
 we are done, so suppose otherwise. Hence, for $i>s$, there is a 
 non-zero homogeneous element $x_i \in R$ such that $x_im_i \in N$. Set $x=x_{s+1} \cdots x_n$.
 Then $xM \subset  N$. By (1), $xM$ is a free graded module. 
 Let $\d=\deg x$. Since $M$ is torsion-free the map $M(\d) \to xM$ given by multiplication by $x$ is 
 an isomorphism. Hence $M(\d)$, and therefore $M$, is a free graded $R$-module.
  \end{pf}

 \begin{cor}
 \label{cor.proj=free}
     Every finitely generated graded projective $C$-module is isomorphic to a direct sum of twists of $C$.
     In particular, a rank one projective graded $C$-module is isomorphic to $C(I)$ for a unique
     $I \in \ZZ_{\rm fin}$.
  \end{cor}
   \begin{pf}
   The only point to be checked is that $C(I) \cong C(J)$ if and only if $I=J$. However, the map
   $$
   \rho:C_{I \oplus J} \to \hom(C(I),C(J)), \qquad \rho(a)(m):= am,
   $$
   is an isomorphism and $\rho(a)$ is an isomorphism if and only if $a$ is a unit, but the only
   homogeneous  units in $C$   are the elements of $k$ which have degree $\varnothing$. Hence $C(I) \cong C(J)$
   if and only if $I \oplus J=\varnothing$, i.e., if and only if $I=J$.
   \end{pf}  
 
 \begin{cor}
 \label{cor.PicX}
   $\Pic \cX \cong \ZZ_{\rm fin}$.
  \end{cor}
   \begin{pf}
  Let $I \in \ZZ_{\rm fin}$. 
  Multiplication in $C$ provides an isomorphism  $C(I) \otimes_C C(I) \stackrel{\sim}{\longrightarrow} C$
  so  $C(I)$ is an invertible object in $\Gr(C,\ZZ_{\rm fin})$. 
 It remains to show that the $C(I)$s are the only invertible objects in $\Gr(C,\ZZ_{\rm fin})$.
 However, if $P \otimes_C Q \cong C$, then $P$ is projective and necessarily of rank one since it 
 embeds in $C$. Hence $P$ is isomorphic to some $C(I)$ by Corollary \ref{cor.proj=free}. 
   \end{pf}

 \begin{prop}
 \label{prop.genors}
 Let $\cS \subset \ZZ_{\rm fin}$. Then the set of projectives $\{Cx_I \; | \; I \in \cS\}$
 generates $\Gr(C,\ZZ_{\rm fin})$ if and only if 
 $$
 \bigcap_{I \in \cS} I = \varnothing
 \quad \hbox{and} \quad
 \bigcup_{I \in \cS} I = \ZZ.
 $$
 \end{prop}
 \begin{pf}
($\Rightarrow$)
By hypothesis there is a surjective map 
$$
\bigoplus_{I \in \cF} Cx_I \to C
$$
for some subset $\cF \subset \cS$. But the image of every non-zero degree preserving homomorphism
 $Cx_I \to C$ is contained in $Cx_I$, so 
$$
\sum_{I \in \cF} Cx_I = C.
$$
Since $C$ is cyclic we can assume $\cF$ is finite. Hence by repeated applications of 
Lemma \ref{lem.K0}(1), the intersection of the $I$s belonging to $\cF$ is empty.

Fix an integer $n$.
By hypothesis there is a surjective map 
$$
\bigoplus_{I \in \cF} Cx_I \to Cx_n
$$
for some subset $\cF \in \cS$. In other words, there is a surjective map 
$$
\bigoplus_{I \in \cF} C(I) \to C(\{n\})
$$
and hence a surjective map 
$$
\bigoplus_{I \in \cF} C(I\oplus\{n\}) \to C.
$$
Since $C(I \oplus \{n\}) \cong Cx_{I \oplus \{n\}}$, it follows from the previous paragraph that the intersection of all the $I \oplus \{n\}$, $I \in \cF$,  is empty.
However, if $n$ does not belong to any of the $I$s  that belong to $\cF$,
$i \in I \oplus \{n\}$ for all $n$, a contradiction.
It follows that $n$ must belong to some $I$. 

($\Leftarrow$)
To prove that $\{Cx_I \; | \; I \in \cS\}$ generates $\Gr(C,\ZZ_{\rm fin})$ it suffices, by Corollary
\ref{cor.proj=free}, to show there is a surjective map
$$
\bigoplus_{I \in \cS} Cx_I \to Cx_K
$$
for every $K \in \ZZ_{\rm fin}$.  By hypothesis, there are elements $I_1,\ldots I_m, I_{m+1},,\ldots, I_n$
of $\cS$ such that $I_1 \cap \cdots \cap I_m \subset K \subset  I_{m+1} \cup \cdots \cup I_n$. 
Write $I= I_1 \cap \cdots \cap I_n$ and $J= I_{1} \cup \cdots \cup I_n$. 
By Proposition \ref{prop.K0}(1), both  $Cx_I$ and $Cx_J$ are quotients of 
$$
\bigoplus_{j=1}^n Cx_{I_j}.
$$
However, $I \cap J \subset K \subset I \cup J$ so Proposition \ref{prop.K0}(2) says that $Cx_K$ is a quotient of $Cx_I \oplus Cx_J$. Hence $Cx_K$ is generated by the $Cx_I$, $I \in \cS$. 
 \end{pf}

 \begin{cor}
  \label{cor.genors1}
 The set of projectives $Cx_n$, $n \in \ZZ$, generates $\Gr C$.
 \end{cor}

  \begin{lem}
 \label{lem.theta1}
 Let $I,J \in \ZZ_{\rm fin}$. Then 
 $$
 \hom_C(C(I),C(J))= C_\varnothing  \theta_{JI}
 $$
 where $\theta_{JI}:C(I) \to C(J)$ is the map $\theta_{JI}(c)=cx_{I \oplus J}$.
 Furthermore, 
  $$
 \bigl(\theta_{JI}\theta_{IJ} - \theta_{JK}\theta_{KJ}\bigr)(c)=\bigl(x_{I \oplus J}^2-x_{J \oplus K}^2\bigr)c.
 $$ 
 \end{lem}
  \begin{pf}
   If $P$ and $Q$ are graded $C$-modules, the $C_\varnothing $-module structure on 
   $\hom(P,Q)$ is given by $(c.f)(p)=f(cp)$ for $c \in C_\varnothing $,   $f \in \hom(P,Q)$, and $p \in P$.  
 It is a standard fact that the map $\rho:C_{I \oplus J} \to \hom_C(C(I),C(J))$ given by
 $$
 \rho(a)(c)=ac
 $$
 is an isomorphism of $C_\varnothing $-modules. Since $C_{I \oplus J}=x_{I \oplus J} C_\varnothing $, 
$ \hom_C(C(I),C(J))$ is generated by $\rho(x_{I \oplus J})$ which is exactly the map $\theta_{JI}$.
The final identity follows immediately from the definition of the $\theta$s.
 \end{pf}

 \subsection{The elements $x_{[n]}$}
 
 We introduce the notation:
\begin{equation}
\label{eq.[n]}
[n]:=  \begin{cases}
	\{1,\ldots,n\} & \text{if $n \ge 1$,}
	\\
	\varnothing &  \text{if $n =0$,}
	\\
	\{{n+1},\ldots, 0\}  & \text{if $n \le -1$.}
	\end{cases}
\end{equation}

\begin{lem}
\label{lem.m+n}
The following identities hold:
 \begin{enumerate}
  \item 
  $[n]=[-n]+n$;
  \item 
  $[n-1]\oplus \{n\}=[n]$;
  \item 
 $[m] \oplus [n]=[n-m]+m$;
 \item{}
 $[m+n]=[m] \oplus ([n]+m)$;
 \item{}
 $-[n]=[-n]-1$.
\end{enumerate}
\end{lem}

  \begin{cor}
  \label{cor.genors2}
 The set  $\{Cx_{[n]} \; | \; n \in \ZZ \}$  generates $\Gr C$.
 \end{cor}
 \begin{pf}
 This follows from the criterion in Proposition \ref{prop.genors}. 
  \end{pf}

 \subsection{Functors between categories of graded modules}
 
 Let $(R,\D)$ and $(S,\G)$ be graded $k$-algebras.
 A {\sf bigraded $R$-$S$-bimodule}  is a $k$-vector space $P$ that is an $R$-$S$-bimodule 
 and has a vector space decomposition
 $$
 P=\bigoplus_{(\d,\c) \in \D \times \G}  P_{(\d,\c)}
 $$
 such that 
 $$
 R_\a .P_{(\d,\c)} .S_\b \subset P_{(\a+\d,\c+\b)}
 $$
 for all $\a,\d \in \D$ and $\b,\c \in \G$.

 For each $\d \in \D$, the subspace 
 $$
 P_{(\d,*)} = \bigoplus_{\c \in  \G}  P_{(\d,\c)}
 $$
 is an $S$-module and we view it as a $\G$-graded $S$-module by declaring that its 
 degree-$\c$ component is $P_{(\d,\c)}$. 
 The left action of an element $r$ in $R_\a$ on $P_{(\d,*)}$ is therefore a degree preserving $S$-module homomorphism $P_{(\d,*)} \to P_{(\a+\d,*)}$, and we therefore obtain a linear map
 \begin{equation}
\label{eq.R->H}
 R_\a \longrightarrow  \hom_S\bigl(P_{(\d,*)},P_{(\a+\d,*)}\bigr).
 \end{equation}
 
 Let $M$ be a $\G$-graded right $S$-module. We define
 $$
 H_S(P,M):= \bigoplus_{\d \in \D} \hom_S(P_{(\d,*)},M)
 $$
 with $\D$-grading given by
  $$
 H_S(P,M)_\d:= \hom_S(P_{(-\d,*)},M).
 $$
Composition of $S$-module homomorphisms gives maps
$$
 \hom_S(P_{(-\d,*)},M) \times  \hom_S\bigr(P_{(-\a-\d,*)},P_{(-\d,*)}\bigr) \to  \hom_S(P_{(-\a-\d,*)},M)
 $$
 and therefore maps
 $$
 H_S(P,M)_{\d} \times R_\a \longrightarrow  H_S(P,M)_{\a+\d}
 $$
 that give $H_S(P,M)$ the structure of a $\D$-graded right $R$-module.

 In summary, we obtain a functor
 $$
 H_S(P,-) : \Gr(S,\G) \to \Gr (R,\D)
 $$
 between the categories of graded right modules.
 A result of del Rio tells us when this is an equivalence of categories. Although 
 del Rio's  result says much more all we need is the following.

 \begin{theorem}
 \cite[Thm. 4.7]{dR} 
 \label{thm.dR}
 With the above notation, suppose that 
 \begin{enumerate}
  \item 
  $P_{(\d,*)}$ is a projective $S$-module for all $\d \in \D$ and
  \item 
  $\{P_{(\d,*)} \; | \; \d \in \D\}$ generates $\Gr(S,\G)$.
\end{enumerate}
Then   $H_S(P,-)$ is an equivalence of categories.
\end{theorem}
  
 \begin{thm}
 \label{thm.main}
 $\Gr A \equiv \Gr(C,\ZZ_{\rm fin})$.
 \end{thm}
 \begin{pf}
 Let $e_nC$, $n \in \ZZ$, be a rank one free $C$-module with basis vector $e_n$ placed in degree $[n]$.
 We define a $\ZZ \times \ZZ_{\rm fin}$ graded vector space $P$ by setting
 $$
 P_{(n,I)} := e_n C_{I \oplus [n]}.
 $$ 
 
Thus $P_{(n,*)} = e_n C$ is  isomorphic to $C([n])$ as a graded right $C$-module.
We give $P$ the structure of an $A$-$C$-module by declaring that  $x$ and $y$ act on $e_nC$ by
$$
x \cdot e_n := e_{n+1} x_{n+1}
\quad \hbox{and}\quad 
y \cdot e_n := e_{n-1} x_{n}.
$$
This {\it does} make $P$ a left $A$-module because 
$$
(xy-yx)e_n= e_n (x_{n}^2-x_{n+1}^2) = e_n .
$$
With this action $P$ is a bigraded $A$-$C$-bimodule. 

  The action of $A_{\ell}$ on $P$ provides a map
 $$
\rho: A_{\ell} \longrightarrow \hom_C(P_{( m,*) }, P_{( \ell +m,*)}).
 $$
Since $P_{(n,*)} \cong C([n])$,  $\hom_C(P_{( m,*) }, P_{( \ell +m,*)})$ is generated as a $C_\varnothing $-module 
by the map $\theta_{[m+\ell],[m]}$ in Lemma \ref{lem.theta1}.

If $\ell \ge 0$, then
$$
x^\ell \cdot e_m  =e_{m+\ell}  x_{m+1} \ldots x_{m+\ell} = e_{\ell+m}  x_{[\ell+m] \oplus [m]}
$$
and
$$
y^\ell \cdot e_m  =e_{m-\ell}  x_{m} \ldots x_{m-\ell+1} = e_{m-\ell}  x_{[m-\ell] \oplus [m]}.
$$
The actions of $x^\ell$ and $y^\ell$ on $P_{(m,*)}$ are therefore the same as the actions of 
$\theta_{[m+\ell],[m]}$ and $\theta_{[m-\ell],[m]}$ respectively. 

Since $xy$ acts on $P_{(n,*)}$ as multiplication by $x_{n}^2$, 
$$
\rho(x^\ell A_0) = \theta_{[m+\ell],[m]} k[x_{n}^2] =  \theta_{[m+\ell],[m]} C_\varnothing .
$$
Similarly, 
$$
\rho(y^\ell A_0) = \theta_{[m-\ell],[m]} k[x_{n}^2] =  \theta_{[m-\ell],[m]} C_\varnothing .
$$
Hence $\rho$ is an isomorphism from $A_\ell$ to $ \hom_C(P_{( m,*)}, P_{(m+ \ell,*)})$.
 
Since this is the case for all $\ell$ and $m$,  and since the $P_{(n,*)}$s provide a set of projective
  generators for  $\Gr(C,\ZZ_{\rm fin})$  
  the theorem follows  from  
  del Rio's  result \cite[Th. 4.7]{dR} (see also \cite[Prop. 2.1]{Sue0}).
  \end{pf}

 \begin{cor}
 $
\Gr A \equiv  \Qcoh \cX.
$
\end{cor}

\section{Simple  and projective graded $C$-modules}
\label{sect.simples}

 As for a Dedekind domain, the classification of all graded $C$-modules follows easily once one has determined the simple and projective ones. By a {\sf simple} graded $C$-module we mean a non-zero graded $C$-module
 whose only graded submodules are itself and the zero module. We define a class of simple graded modules that we call {\it special}. These are the simple $\cO_\cX$-modules that are supported at the stacky points
 of $\cX$. The importance of these modules is apparent from Proposition \ref{prop.proj} and 
 Corollary \ref{cor.simples.funct}.

Under the equivalence $\Qcoh \cX \equiv  \Gr(C,\ZZ_{\rm fin})$ the projective graded $C$-modules
correspond to  the locally free  $\cO_\cX$-modules.

\subsection{The simple graded modules}

\begin{prop}
The maximal graded ideals of $C$ are the ideals $(p)$ as $p$ ranges over the irreducible 
elements in $C_\varnothing -\{x_n^2 \; | \; n \in \ZZ\}$ and the ideals $(x_n)$ for $n \in \ZZ$.
\end{prop}
\begin{pf}
The ideals in the statement of the proposition are certainly graded ideals. 

To show a graded ideal $\fa$ is maximal among graded ideals it suffices to show that every homogeneous element of $C-\fa$ is a unit in $C/\fa$.
Every homogeneous element of $C$ is of the form $fx_I$ for some $f \in C_\varnothing $ and 
some  $I \in \ZZ_{\rm fin}$. 

Let $\fa = (p)$ where $p$ is an irreducible element of $C_\varnothing $ but not one of the $x_n^2$. 
Every element of $C_\varnothing  - \fa$ is a unit modulo $\fa$. 
If $i \in \ZZ$, then $x_i^2 \in C_\varnothing  -(p)$ so $x_i^2$,
and therefore $x_i$, is a unit modulo $\fa$. It follows that every $x_I$ is a unit modulo $\fa$. 
Hence if $f \in C_\varnothing $ and $fx_I \notin \fa$, then $fx_I$ is a unit modulo $\fa$.
This completes the proof that $(p)$ is a maximal graded ideal.

Now let $\fa=(x_n)$. If $i \in \ZZ-\{n\}$, then $x_i^2$ is congruent to a non-zero scalar modulo $\fa$ so is 
a unit. Since $C_\varnothing =k[x_n^2]$, $\fa \cap C_\varnothing $ is a maximal ideal of $C_\varnothing $. 
Hence every element in $C_\varnothing -\fa$ is a unit modulo $\fa$. Therefore, if $f \in C_\varnothing $ and 
$fx_I \notin \fa$, then
$fx_I$ is a unit modulo $\fa$. This completes the proof that $(x_n)$ is a maximal graded ideal for all
 $n \in \ZZ$.

Now let $\fa$ be an arbitrary maximal graded ideal of $C$. Then $(C/\fa)_\varnothing$ is a field, so $\fa$ contains
an irreducible element of $C_\varnothing $, say $p$, and $\fa \supset (p)$. If $p \not\in \{x_n^2 \; | \; n \in \ZZ\}$, then
$(p)$ is maximal so $\fa=(p)$. On the other hand, if $p=x_n^2$, then $x_n \in \fa$ because $C/\fa$ has no homogeneous zero divisors, and therefore $\fa=(x_n)$.  
\end{pf}

\subsection{The ordinary simple graded $C$-modules}

The simple graded modules of the form $C/(p)$ where $p$ is an irreducible element of $C_\varnothing $
play virtually no role in this paper. However, the following facts are easily verified:
\begin{enumerate}
  \item 
 if $p$ and $p'$ are relatively prime  irreducibles, then $C/(p) \not\cong C/(p')$
 and $\ext^1_C(C/(p),C/(p')) =0$;
  \item 
  for all $J \in \ZZ_{\rm fin}$, $(C/(p))(J) \cong C/(p)$;
 \item{}
 $\ext^1_C(C/(p),C/(p)) \cong C/(p)$. 
\end{enumerate}
The simple modules of the form $C/(p)$ correspond to the non-stacky points of the coarse moduli space for 
$\cX$. 

For each $\l \in k-\ZZ$,  we define
$$
\cO_\l:={{C}\over{(x_0^2-\l)}}.
$$

\subsection{The special simple graded $C$-modules}
\label{sect.specials}
A simple graded $C$-module is {\sf special} if it is isomorphic to one of the modules
$$
X_n := {{C}\over{(x_n)}}, 
 \quad Y_n  :=  \biggl({{C}\over{(x_n)}}\biggr) \bigl(\{n\}\bigr), \quad n \in \ZZ.
 $$
The following observations follow immediately from the definition:
\begin{enumerate}
\item
There are non-split exact sequences 
$$
\UseComputerModernTips
  \xymatrix{
  0 \ar[r]  & Y_n  \ar[r] & {{C}\over{(x_n^2)}} \ar[r] & X_n   \ar[r] & 0.
  }
  $$
  and 
$$
\UseComputerModernTips
  \xymatrix{
  0 \ar[r]  & X_n  \ar[r] &\biggl( {{C}\over{(x_n^2)}}\biggr)(\{n\}) \ar[r] & Y_n   \ar[r] & 0
  }
  $$
\item
As $C_\varnothing $-modules, the homogeneous components of $X_n$ are
$$
(X_n)_I \cong \begin{cases}
	C_\varnothing /(x_n^2) & \text{if $n \notin I$,}
	\\
	0 & \text{if $n \in I$.}
	\end{cases}
$$
\item
As $C_\varnothing $-modules, the homogeneous components of $Y_n$ are
$$
(Y_n)_I \cong \begin{cases}
	C_\varnothing /(x_n^2) & \text{if $n \in I$,}
	\\
	0 & \text{if $n \notin I$.}
\end{cases}
$$
  \item 
  $Y_n \not\cong X_n$ because $(X_n)_\varnothing \cong k$ but $(Y_n)_\varnothing=0$, and 
  \item 
  $Y_n(\{n\}) =X_n$ because $2\{n\}=0$.
  \end{enumerate}

One may define and/or characterize the special simple modules in terms of their properties
inside the category $\Gr(C,\ZZ_{\rm fin})$. For example, working with $\Gr A$, Sierra characterizes them as the  simple graded modules $S$ for which $\ext^1(S,M) \ne 0$ for some simple graded module $M \not\cong S$. In order
to offer an alternative to Sierra's characterization we will characterize them as those simples $S$ for which
$\hom(P,S)=0$ for some non-zero projective graded module $P$ (Proposition \ref{prop.specials}).

As we shall see, the isomorphism class of a special simple module is determined by the degrees
in which it is zero, and a simple graded module is special if and only if some of its homogeneous 
components are zero.

\begin{prop}
\label{prop.simples}
Let $I \in \ZZ_{\rm fin}$ and let $n \in \ZZ$.  
\begin{enumerate}
  \item 
  $\hom(C,X_n) \cong k$.
\item{}
 $ \hom(C,Y_n)=0$.
  \item 
$
X_n(I) \cong
\begin{cases}
X_n & \text{if $n \notin I$,}
\\
Y_n & \text{if $n \in I$.}
\end{cases}
$
  \item 
  $
\hom(Cx_I,X_n)  \cong 
	\begin{cases}
k & \text{if $n \notin I$,}
\\
0 & \text{if $n \in I$.}
\end{cases}
$
\item
$I=\{n \in \ZZ \; | \; \hom(Cx_I,X_n) =0\}$.
\item
$I= \{n \in \ZZ \; | \; \hom(Cx_I,Y_n) \ne 0\}$.
\end{enumerate}
\end{prop}
\begin{pf}
(1)
We have $\hom(C,X_n) \cong (X_n)_\varnothing \cong  C_\varnothing /C_\varnothing  x_n^2 \cong k$.

(2)
We have $\hom(C,Y_n) \cong (Y_n)_\varnothing =(X_n)_{\{n\}} =0$. 

(3)
If $m \ne n$, then the image of $x_m$ in $C/(x_n)$ is a unit so  multiplication by $x_m$ 
is an isomorphism $X_n(\{m\}) \cong X_n$. In general, if $I=\{i_1,\ldots,i_t\}$, then 
$$
X_n(I)=X_n(\{i_1\}) \cdots (\{i_t\})
$$
so the result follows from the previous sentence. 

(4)
Since $\hom(Cx_I,X_n) \cong \hom(C,X_n(I))$,  this follows from (1) and (3).

(5) and (6) follow from (4).
\end{pf}

\begin{prop}
\label{prop.specials}
Let $S$ be a simple graded module. The following three conditions are equivalent:
\begin{enumerate}
  \item 
  $S$ is special;
  \item 
  $\hom(P,S) =0$ for some non-zero projective graded module $P$;
  \item 
  $S_J=0$ for some $J \in \ZZ_{\rm fin}$.
\end{enumerate}
\end{prop}
\begin{pf}
By Corollary \ref{cor.proj=free},  every projective graded $C$-module is a direct sum of various $C(J)$s, so 
(2) holds  if and only if $\hom(C(J),S) =0$ for some $J \in \ZZ_{\rm fin}$. 
However, $\hom(C(J),S) \cong S_J$ so (2) holds  if and only if $S_J=0$ for some $J \in \ZZ_{\rm fin}$. 
This proves the equivalence of (2) and (3).

Suppose $S$ is not special. Then $S\cong C/pC$ for some irreducible $p \in C_\varnothing $ and $S_J = C_J/pC_J$ 
for all $J \in \ZZ_{\rm fin}$. But $C_J$ is isomorphic to $C_\varnothing $ as a $C_\varnothing $-module so $S_J \ne 0$
for all $J \in \ZZ_{\rm fin}$. On the other hand, if $S$ is special, then $S_J$ {\it is} zero for some $J$ by
parts (5) and (6) of Proposition \ref{prop.simples}. This proves the equivalence of (1) and (3).
\end{pf}

The next result corresponds to Sierra's result \cite[Thm. 5.5]{Sue1}. Our proof is a little different. For example, we characterize the special simple graded modules $S$ using Proposition
\ref{prop.specials}(2), and we also exploit the fact that $C$ is commutative by using the map 
$\Pic(C,\ZZ_{\rm fin}) \to \Aut(C_\varnothing )$, $[F] \mapsto \a_F$, defined in Proposition \ref{prop.aut.R0}.

We write $\Iso(\ZZ)$ for the isometry group of $\ZZ$ with respect to  the metric is $d(m,n)=|m-n|$. 
The isometries are exactly the maps $n \mapsto \ve n+d$ where 
$\ve =\pm 1$ and $d \in \ZZ$.  
As an abstract group, $\Iso(\ZZ)$ is isomorphic to the dihedral group  $D_\infty$.

\begin{thm}
\label{thm.permute}
There is a  group homomorphism
\begin{equation}
\label{eq.Pic.Iso}
\Pic(C,\ZZ_{\rm fin}) \to \Iso(\ZZ), \qquad [F] \mapsto (n \mapsto \ve n+d),
\end{equation}
where $\ve \in \{ \pm 1\}$ and $d \in \ZZ$ are determined by the requirement that 
\begin{equation}
\label{eq.F.isometry}
FX_n \oplus FY_n \cong X_{\ve n+d} \oplus Y_{\ve n+d}
\end{equation}
for all $n \in \ZZ$  and 
\begin{equation}
F\cO_\l \cong \cO_{\ve \l +d}
\end{equation}
for all $\l \in k-\ZZ$.  
\end{thm}
\begin{pf}
Let $S$ and $S'$ be special simple graded modules. By parts (5) and (6) of Proposition \ref{prop.simples},
 $S \oplus S'  \cong X_m \oplus Y_m$  for some $m \in \ZZ$ if and only if 
 $\dim_k\hom(C(I),S \oplus S')  =1$ for all $I \in \ZZ_{\rm fin}$. Since 
 $$
 \dim_k  \hom(F(C(I)),FX_n \oplus FY_n) = \dim_k\hom(C(I),X_n \oplus Y_n) 
 $$
 and since $F$ permutes the isomorphism classes of rank one projective graded modules, it follows that
 $FX_n \oplus FY_n \cong   X_{g(n)} \oplus Y_{g(n)}$ for a unique $g(n) \in \ZZ$.
 Since $F$ is an auto-equivalence $g$ is a permutation of $\ZZ$.   
 
 By Proposition \ref{prop.aut.R0}, $F$ determines an automorphism $\a_F$ of $C_\varnothing $ 
 having the property that $a \in C_\varnothing $ annihilates a module $M$ if and only if $\a_F(a)$ annihilates
 $FM$.  Since $x_n^2$ annihilates $X_n \oplus Y_n$ , $\a_F(x_n^2)$ is a multiple of $x_{g(n)}^2$.

Write $z=x_0^2$.  Since $C_\varnothing $ is the polynomial ring $k[z]$, there is $\ve \in k-\{0\}$ and $d \in k$ such that
$\a_F(z)=\ve z +d$.  Therefore $\ve z -n+d$, which is $\a_F(z-n)$, is a scalar multiple of 
$z-g(n)$. Hence $\ve z-n+d=\ve(z-g(n))$
for all $n \in \ZZ$.  Thus $g(n)={{1}\over{\ve}}(n-d)$ for all $n \in \ZZ$. It follows that $\ve=\pm 1$
 and $d \in \ZZ$.   
\end{pf}

Theorem \ref{thm.Pic}  shows that the map (\ref{eq.F.isometry}) is surjective and that its kernel 
is the image of $\ZZ_{\rm fin}$.

\subsection{Projective graded modules}
 
By Corollary \ref{cor.proj=free},  every projective graded $C$-module is a direct sum of various $C(I)$s.
The next two results are immediate consequences of parts (5) and (6) of Proposition \ref{prop.simples}.

 \begin{cor}
  If $P$ is a rank one projective graded $C$-module, then
  $$
  \hom(P,X_n) \ne 0 \; \Longleftrightarrow \;  \hom(P,Y_n)=0.
  $$
   \end{cor}

  \begin{cor}
 \label{cor.projs+simples2}
 A rank one projective graded $C$-module maps surjectively onto infinitely many $X_n$s but
 only finitely many $Y_n$s.
  \end{cor}

  \begin{remark}
 By Corollary \ref{cor.projs+simples2}, the $X_n$s play a different role in 
  $\Gr(C,\ZZ_{\rm fin})$ from the $Y_n$s. The $X_n$'s are the simple $G$-equivariant 
  $C$-modules on which the corresponding isotropy groups act trivially, 
  whereas those isotropy groups act on the $Y_n$s via the sign representation.
  
   Sierra labels the simple graded $A$-modules $X(n)$ and $Y(n)$, $n \in \ZZ$, but her labelling is {\it not}
  compatible with ours---her $X(n)$ corresponds to $X_n$ if $n \ge 0$ and to $Y_n$ if $n<0$. 
     Her labelling, which is designed to remind the reader that $X(n)$ (resp., $Y(n)$) 
   is isomorphic as an ungraded $A$-module to $A/xA$ (resp., $A/yA$)  makes the different properties of 
    the $X_n$s and $Y_n$s less apparent.  
   \end{remark}

\begin{prop}
\label{prop.proj}
Let $P$ and $Q$ be finitely generated projective graded $C$-modules having the same rank. 
Then $P \cong Q$ if and only if  $\dim \hom_C(P,Y_n) =  \dim \hom_C(Q,Y_n)$  for all $n$.
\end{prop}
\begin{pf}
($\Leftarrow$)
Suppose $\dim \hom_C(P,Y_n) =  \dim \hom_C(Q,Y_n)$  for all $n$. We will argue by induction on  
$r:=\rank P$. 

By Corollary \ref{cor.proj=free} and  parts (5) and (6) of Proposition \ref{prop.simples},
the result is true when $r=1$, so we assume that $r \ge 2$.

Let $\fa$ be the largest graded ideal of $C$ that is the image of a degree-preserving homomorphism
$P \to C$.  Since $P \cong Q$, $\fa$ is also the largest graded ideal of $C$ that 
is the image of a degree-preserving homomorphism $Q \to C$. 
Since $\fa$ is projective, there are graded projectives $P'$ and $Q'$ of the same rank such that 
$P \cong P' \oplus \fa$ and $Q \cong Q' \oplus \fa$. 
It is obvious that $\dim \hom_C(P',Y_n) =  \dim \hom_C(Q',Y_n)$  for all $n$ so, by the induction hypothesis, 
$P' \cong Q'$. It follows that $P \cong Q$. 

($\Rightarrow$)
This is obvious. 
\end{pf}

\begin{cor}
\label{cor.simples.funct}
Let $F$ and $G$ be autoequivalences of $\Gr(C,\ZZ_{\rm fin})$. Then $F \cong G$ if and 
only if $FS \cong GS$
for all special simples $S$.
\end{cor}
\begin{pf}
Suppose  $FS \cong GS$ for all special simples $S$. 
Then for every finitely generated graded $C$-module
$M$ and every special simple $S$,
\begin{align*}
\dim \hom(FM,FS) & =\dim \hom(M,S) 
\\
& =\dim \hom(GM,GS)
\\
&= \dim \hom(GM,FS).
\end{align*}
But $F$ and $G$ permute the special simples by Theorem \ref{thm.permute}, so 
$$
\dim \hom(FM,X_i)= \dim  \hom(GM,X_i)
$$
and   $$
\dim \hom(FM,Y_i)= \dim \hom(GM,Y_i)
$$
for all $i$.  Now take $M=R(j)$, $j \in \G$. Because $F(R(j))$ and $G(R(j))$ are rank one projectives,  
it follows from Proposition  \ref{prop.proj} that $F(R(j)) \cong G(R(j))$.   
Proposition \ref{prop.wisom} now implies  that $F \cong G$.
 \end{pf}

\begin{cor}
\label{cor.FG}
Let $\sA$ be a  category and suppose that $F,G:\Gr(C,\ZZ_{\rm fin}) \to \sA$ are equivalences of categories. 
Then $F \cong G$ if and only if $FS \cong GS$ for all special simples $S$.
\end{cor}
\begin{pf}
Let $G^{-1}$ be a quasi-inverse to $G$. Then $F \cong G$ if and only if $G^{-1}F \cong \id_{\Gr C}$.
Hence $F \cong G$ if and only if $G^{-1}FS \cong S$ for all special simples $S$. The result follows.
 \end{pf}

\section{The Grothendieck group of $\cX$}
\label{sect.K0}

The Grothendieck group of $\cX$, denoted $K_0(\cX)$, is, by definition,  the Grothendieck 
group of the category of locally free coherent $\cO_\cX$-modules. 
Under the equivalence $\Qcoh \cX \equiv  \Gr(C,\ZZ_{\rm fin})$ locally free  coherent
$\cO_\cX$-modules correspond to finitely generated projective graded $C$-modules.  

We write $\gr(C,\ZZ_{\rm fin})$, or $\gr C$, 
and $\sP$ respectively for the full subcategories of $\Gr(C,\ZZ_{\rm fin})$ consisting of the  finitely generated modules and the finitely generated projective graded modules.

Because $C$ is a graded principal ideal domain, every graded $C$-module $M$ is isomorphic to 
$P/Q$ where $P$ and $Q$ are projective, even free, graded modules.  The natural map 
$K_0(\sP) \to K_0(\gr C)$ is therefore an isomorphism.\footnote{Hence the natural map 
$K_0(\cX) \to K_0(\coh \cX)$ is also an isomorphism.}
In particular, 
$$
K_0(\cX) \cong K_0(\gr C).
$$

\subsection{Classification of projective graded modules}

If $I$ and $J$ are multi-sets, i.e., sets whose elements have multiplicities, their union as a multi-set will
be denoted by $I \boxplus J$.

\begin{prop}
\label{prop.multi}
Let $I_1,\ldots,I_r,J_1,\ldots,J_r \in \ZZ_{\rm fin}$. Then 
$$
\bigoplus_{n=1}^r C(I_n) \cong  \bigoplus_{n=1}^r C(J_n)
$$
if and only if
$$
I_1 \boxplus \cdots \boxplus I_r =  J_1 \boxplus \cdots \boxplus J_r.
$$
\end{prop}
\begin{pf}
By Proposition \ref{prop.simples}(6), $\dim\hom(C(J),Y_i)$ is 1 if $i \in J$ and $0$ otherwise, so the result follows from Proposition \ref{prop.proj}.
\end{pf}

Let $r \in \NN$ and write $\ZZ_{{\rm mult}\le r}$ for the set of all finite multi-sets $M$ 
of integers such that every  element of $M$ has multiplicity $\le r$. Define
$$
\Phi_r: \ZZ_{{\rm mult}\le r} \to \Gr(C,\ZZ_{\rm fin})
$$
by declaring that
$$
\Phi_r(M):= C(I_1 )\oplus \cdots \oplus C(I_d) \oplus C^{r-d}
$$
where $I_1,\ldots,I_d$ are the unique elements of $\ZZ_{\rm fin}$ such that 
$$
I_1 \supset \cdots \supset I_d \ne \varnothing
\quad \hbox{ and} \quad I_1 \boxplus \cdots \boxplus I_d = M.
$$

\begin{cor}
Fix a non-negative integer $r$. Then $\Phi_r$ gives a bijection between the elements of 
$\ZZ_{{\rm mult}\le r}$  and  the isomorphism classes of  finitely generated projective graded 
$C$-modules of rank $\le r$. The inverse to $\Phi_r$ sends a module isomorphic to 
$C(J_1 )\oplus \cdots \oplus C(J_r)$ to $ J_1 \boxplus \cdots \boxplus J_r$.
\end{cor}

If $P$ is a finitely generated projective graded $C$-module we write $[P]$ for its class in the 
Grothendieck group $K_0(\gr C)$.

  \begin{cor}
 \label{cor.K0}
 Let $P$ and $Q$ be finitely generated projective graded $C$-modules of the same rank. 
 Then $[P]=[Q]$  if and only if $P \cong Q$.   
  \end{cor}
  \begin{pf}
If $[P]=[Q]$, there is a finitely generated graded projective $M$ such that $P \oplus M \cong Q \oplus M$.
It follows that $\dim \hom_C(P,Y_n) =  \dim \hom_C(Q,Y_n)$  for all $n$ so $P \cong Q$.
The reverse implication is trivial.
  \end{pf}

 Because $C$ is commutative, $K_0(\gr C)$ is a commutative ring with product $[P].[Q]=[P\otimes_C Q]$
 where the tensor product is the usual tensor product of graded modules. 
  By Proposition \ref{prop.K0},  
$$
[C(\{m\})].[C(\{n\})] =
\begin{cases}
[C] & \text{if $m=n$}
\\
 [C(\{m\})]  + [C(\{n\})]  -[C] & \text{if $m \ne n$.}
 \end{cases}
$$
 In due course, we will see that the classes $[\CC(\{n\})]$, $n \in \ZZ$, provide a $\ZZ$-basis for 
 $K_0(\gr C)$.

 \subsection{The homomorphism $\Upsilon$}  
  As in section \ref{sect.gp.sch}, we write $u_I$ for the element of   the integral group ring
  $\ZZ\ZZ_{\rm fin}$ corresponding to $I$.  There is a surjective  ring homomorphism
 $$
\Upsilon:  \ZZ\ZZ_{\rm fin} \to K_0(\gr C), \quad u_I \mapsto [C(I)].
 $$

 \begin{thm}
 \label{thm.K0.relns}
 The kernel of the homomorphism $\Upsilon : \ZZ\ZZ_{\rm fin} \to K_0(\gr C)$ is the ideal  generated by  
 the elements
 $$
 u_I + u_J - u_{I \cap J} - u_{I \cup J}, \qquad I,J \in \ZZ_{\rm fin}.
 $$
 Equivalently, $\ker(\Upsilon)$ is generated by    $\{u_mu_n+1 -u_m - u_n \; | \; m \ne n\}$.
 \end{thm}
 \begin{pf}
 Let $\fa$ be the ideal generated by the elements  $u_I + u_J - u_{I \cap J} - u_{I \cup J}$.
 Then $\fa \subset \ker \Upsilon$ because 
 \begin{equation}
 \label{eq.I+J}
 C(I) \oplus C(J) \cong C(I \cap J) -C(I \cup J).
 \end{equation}
 
 In order to shorten the notation we will write $I+J$ rather than $u_I + u_J$ in this proof. 
 
Let $x \in \ker \Upsilon$. By Proposition \ref{prop.multi}  and Corollary \ref{cor.K0}, 
$$
x=(I_1+ \cdots + I_n) -(J_1+\cdots +J_n)
$$
for some $n\in \NN$ and some elements $I_r, J_r \in \ZZ_{\rm fin}$ with the property that 
$$
I_1 \boxplus \cdots \boxplus I_n =  J_1 \boxplus \cdots \boxplus J_n.
$$
It follows that $I_1 \cap \ldots \cap I_n=J_1 \cap \ldots \cap J_n$ and 
 $I_1 \cup \ldots \cup I_n=J_1 \cup \ldots \cup J_n$.

We will argue by induction on $n$ to show that $x \in \fa$. If $n \le 1$, there is nothing to prove and when
 $n=2$ the result follows from  (\ref{eq.I+J}). Suppose that $n \ge 3$. 
 
A sequence of elements $K_1,\ldots,K_r$ belonging to $\ZZ_{\rm fin}$ is said to be {\sf decreasing} 
if $K_1 \supset K_2 \supset \ldots \supset K_r$.  
There is a unique decreasing sequence $K_1,\ldots,K_n$ such that
$$
I_1 \boxplus \cdots \boxplus I_n =  K_1 \boxplus \cdots \boxplus K_n.
$$
By Proposition \ref{prop.multi}, $(I_1+ \cdots + I_n) -(K_1+\cdots +K_n)$ and 
$(J_1+ \cdots + J_n) -(K_1+\cdots +K_n)$ belong to $\ker \Upsilon$. 
If   $(I_1+ \cdots + I_n) -(K_1+\cdots +K_n)$ and 
$(J_1+ \cdots + J_n) -(K_1+\cdots +K_n)$ belong to $\fa$ so does $x$.
It therefore suffices to show that $x$ belongs to $\fa$ when $I_1 ,\dots,I_n$ is decreasing. 
We assume that is the case.

Let $L:=I_1 \cap \ldots \cap I_n$ and define $I_s':=I_s-L$ and $J_s':=J_s - L$ for $1 \le s \le n$, and 
write 
$$
x'=(I'_1+ \cdots + I'_n) -(J'_1+\cdots +J_n').
$$ 
Notice that $x=u_Lx'$.

It is clear that  $I'_1, \ldots,I'_n$ is decreasing, $I'_1 \supset J_1'$, $I_n'=\varnothing$, and 
$$
I_1' \boxplus \cdots \boxplus I_n' =  J_1' \boxplus \cdots \boxplus J_n'.
$$
Let $K=I_1' - J_1'$. Then 
$$
K \boxplus I_2' \boxplus \cdots \boxplus I_{n-1}'=  J_2' \boxplus \cdots \boxplus J_n'
$$
so 
$$
y:=(K+I_2' \cdots + I_{n-1}')   -(J_2' + \cdots +J_n')
$$
belongs to $\ker \Upsilon$. It follows from the  induction hypothesis that $y \in \fa$.
But $(I_1'+\varnothing) - (J_1'+L) \in \fa$ and $x'=y+(I_1'+\varnothing) - (J_1'+L)$, so $x' \in \fa$ too. 
Since $x=x'u_L$, $x \in \fa$.
 \end{pf}

\begin{cor}
\label{cor.K0.basis}
The elements $[C]$ and $[C(\{m\})]$, $m \in \ZZ$, provide a $\ZZ$-basis for $K_0(\gr C)$.
\end{cor}
\begin{pf}
Let $\fa=\ker \bigl(\Upsilon:\ZZ\ZZ_{\rm fin} \to K_0(\gr C)\bigr)$.  
Since $u_mu_n \equiv  u_m+u_n-1$ modulo $\fa$, it follows that 
$\ZZ\ZZ_{\rm fin}/\fa$ is spanned by the images of $u_m$, $m \in \ZZ$, and 1.  

If there is a relation in $K_0(\gr C)$ of the form
$$
[C(\{m_1\})]+ \cdots + [C(\{m_r\})]+ d [C]  = [C(\{n_1\})]+ \cdots + [C(\{n_s\})]+ e [C] 
$$
for some positive integers $d$ and $e$, and elements $m_i$ and $n_j$ in $\ZZ$, then 
$r+d=s+e$ and the multisets $\{\{m_1,\ldots,m_r\}\}$ and  $\{\{n_1,\ldots,n_s\}\}$
are equal. Hence $d=e$, and it follows from this that the images of $u_m$, $m \in \ZZ$, and 1
in $\ZZ\ZZ_{\rm fin}/\fa$ are linearly independent.  
\end{pf}

  \begin{cor}
 If $\l \in k-\ZZ$, then $[\cO_\l]=0$ but $[X_n] = -[Y_n] \ne 0$ for every $n \in \ZZ$.
 \end{cor}
\begin{pf}
Since $\cO_\l = C/(x_0^2-\l)$ and $\deg( x_0^2-\l)=\varnothing$, $[\cO_\l]=0$. On the other hand, since the $[C(\{n\})]$s form a basis and $X_n = C/Cx_n$, $[X_n]=[C]-[C(\{n\})] \ne 0$. By definition, $Y_n=X_n(\{n\})$, so $[Y_n]=-[X_n]$.
\end{pf}

\section{Symmetries and automorphisms of $\Gr(C,\ZZ_{\rm fin})$}
\label{sect.symm}
\label{sect.Omega}

Consider the diagram
$$
\UseComputerModernTips
  \xymatrix{
\cdots &\ar@{-}[l] \!\! \colon \! \! \! \!  \ar@{-}[r] &  \!\! \colon \! \! \! \!  \ar@{-}[r] & \!\! \colon \! \! \! \!  \ar@{-}[r] & \!\! \colon \! \! \! \!  \ar@{-}[r]  &\!\! \colon \! \! \! \! \ar@{-}[r]& \cdots
 }
  $$ 
in which the underlying line is $\Spec C_\varnothing =\Spec k[z]$ and the two fractional points at the loci $x_n^2=0$, 
$n \in \ZZ$, represent the special simples $X_n$ and $Y_n$. 
There are two obvious symmetries: translation $n \mapsto n+1$, and reflection about $0$. 
The automorphism $z \mapsto z+1$ of $k[z]$ extends to an automorphism $\tau$ of $C$, and 
the automorphism $z \mapsto -z$ of $k[z]$ extends to an almost-automorphism $\varphi$ of $C$.
(If $\sqrt{-1} \in k$, the automorphism $z \mapsto -z$ of $k[z]$ extends to an automorphism 
$x_n=\omega:\sqrt{z-n} \mapsto \sqrt{z+n}=\sqrt{-1}x_{-n}$ of $C$ such that $\omega_*=\varphi_*$.)

\begin{thm}
\label{thm.symm}
 \label{thm.tau*.simples} 
There is an automorphism $\tau$  and an almost-automorphism $\varphi$ such that
$$
\begin{matrix}
\tau_* X_n \cong X_{n+1}   & \; &   \tau_*Y_n \cong Y_{n+1}   & \; &   \tau_*\cO_\l \cong \cO_{\l+1}
\\
&
\\
\varphi_*X_n \cong X_{-n}  &\; & \varphi_* Y_n \cong Y_{-n}   & \; &   \varphi_*\cO_\l \cong \cO_{-\l}.
\end{matrix}
$$
for all $n \in \ZZ$ and $\l \in k-\ZZ$.
\end{thm}
\begin{pf}
Define the automorphism  $\tau$  by 
$$
\tau(x_n)=x_{n-1}, \qquad n \in \ZZ.
$$
Let $\tau_*$ be the associated automorphism  of $\Gr(C,\ZZ_{\rm fin})$
defined in section \ref{sect.autom.twist}. 

Then $\tau_* X_n$ is simple, and its degree $I$ component, $(\tau_* X_n)_I$, is equal
to $(X_n)_{\tauol I}$ which is $(X_n)_{I-1}$. By section \ref{sect.specials}, 
$(X_n)_{I-1}$  is zero exactly 
when $n \in I-1$, i.e., when $n+1 \in I$. Hence $\tau_* X_n \cong X_{n+1}$. 
The proof for $\tau_* Y_n$ is similar. 
Since $\cO_\l$ is annihilated by $x_0^2-\l$, $\tau_*\cO_\l$ is annihilated by $\tau^{-1}(x^2_0-\l) = x_{1}^2 -\l
= x_0^2-1-\l$. Hence $\tau_*\cO_\l \cong \cO_{\l+1}$.
 
The existence of $\varphi$ is proved in Proposition \ref{prop.varphi} below. 
There is an associated automorphism  $\varphi_*$ of $\Gr(C,\ZZ_{\rm fin})$.  
Only two properties of $\varphi$ 
are needed for the proof this theorem: $\varphi(C_I)=C_{-I}$ for all $I \in \ZZ_{\rm fin}$ 
and $\varphi(x_0^2)=-x_0^2$. 
Because $\varphi(C_I)=C_{-I}$ the degree $I$ component  of $\varphi_* X_n$, which is, 
of course, a simple graded $C$-module, is equal   to $(X_n)_{-I}$. 
By section \ref{sect.specials}, $(X_n)_{-I}$ is zero exactly 
  when $n \in -I$, i.e., when $-n \in I$. Hence $\varphi_* X_n \cong X_{-n}$. 
  The proof for $\varphi_* Y_n$ is similar. Finally, since  $\cO_\l$ is annihilated by $x_0^2-\l$, 
  $\varphi_*\cO_\l$ is annihilated by $\varphi^{-1}(x_0^2-\l) = -x_0^2-\l$.
\end{pf}

\begin{prop}
\label{prop.varphi}
Write $z=x_0^2$. There is an almost-automorphism $\varphi:C \to C$ defined by the conditions:
\begin{itemize}
  \item 
   $\varphi:k[z] \to k[z]$ is the $k$-algebra automorphism  $\varphi(z)=-z$ and
  \item 
  $
\varphi(ax_J)=\varphi(a)x_{-J}
$
for all $a \in k[z]$ and $J \in \ZZ_{\rm fin}$. 
\end{itemize}
Furthermore, 
\begin{enumerate}
  \item 
   if $c \in C_I$ and $d \in C_J$, then  
$
\varphi(cd)=(-1)^{|I \cap J|} \varphi(c) \varphi(d);
$
  \item 
  $
\varphi(x_n^2)= - x_{-n}^2
$
for all $n \in \ZZ$;
\item
$\varphi^2=\id_C$.
\end{enumerate}
\end{prop}
\begin{pf}
Let $I,J,K \in \ZZ_{\rm fin}$. We write $\l_{I,J}:= (-1)^{|I \cap J|}$. 
Since 
$$
|K \cap (I \oplus J)| =|K \cap I| + | K \cap J| \; \hbox{(mod 2)}
$$
it follows that 
$$
|K \cap (I \oplus J)| + |I \cap J| =  | (K \oplus I) \cap J|  + |K \cap I|  \; \hbox{(mod 2)}
$$
and hence that
$$
\l_{K,I \oplus J} \l_{I,J} = \l_{K \oplus I,J}\l_{K,I}.
$$
To show that $(\varphi,\l)$ is an almost-automorphism of $C$  it therefore suffices to prove (1).
But first we observe that (2) is true because  $x_n^2=z-n$. 

(1)
It is enough to check this for $c=x_I$ and $d=x_J$. In that case
\begin{align*}
\varphi(x_I x_J) & = \varphi(x_{I \cap J}^2 x_{I \oplus J})
\\
 & = (-1)^{|I \cap J|} x^2_{-(I \cap J)}  x_{-(I \oplus J)}
 \\
 & = (-1)^{|I \cap J|}  x_{-I} x_{-J}
 \\
 & = (-1)^{|I \cap J|}  \varphi(x_I)\varphi( x_J).
 \end{align*}

 (3)
 This is clear.
\end{pf}

\begin{prop}
\label{prop.omega*}
If $\sqrt{-1} \in k$, there is an algebra automorphism 
$$
\omega: C \to C, \qquad  \omega(x_n):=\sqrt{-1} x_{-n},
$$
such that
$$
\omega_*  \cong \varphi_*.  
$$
\end{prop}
\begin{pf}
It is easy to see that $\omega$ does extend to an algebra automorphism. 
To show that $\omega_* \cong \varphi_*$ it suffices to show that their actions on isomorphism
classes of the special simple graded modules are the same.   However, $\omegaol = \varphiol$ 
because $\omega(C_I)=C_{-I} = \varphi(C_I)$ for all $I \in \ZZ_{\rm fin}$ so 
the same argument as was used in Theorem \ref{thm.tau*.simples} for the action of
 $\varphi_*$ on the special simples shows that 
$\omega_* X_n \cong X_{-n}$ and $\omega_* Y_n \cong Y_{-n}$.
The result follows.
\end{pf}

\begin{thm}
\label{thm.Pic}
There is an exact sequence
$$
1 \to \ZZ_{\rm fin}  \to \Pic(C,\ZZ_{\rm fin}) \to \Iso(\ZZ) \to 1
$$
where the map into $\Pic(C,\ZZ_{\rm fin})$ sends $J$ to the twist functor $(J)$ and the map
out of $\Pic(C,\ZZ_{\rm fin})$ is described in Theorem \ref{thm.permute} (equivalently, it 
sends $F$ to the automorphism  $\a_F$ of $C_\varnothing $ defined in Proposition \ref{prop.aut.R0}).
\end{thm}
\begin{pf}
Since $\Iso(\ZZ)$ is generated by the maps $n \mapsto n+1$ and $n \mapsto -n$, it follows from Theorem 
\ref{thm.tau*.simples} that the map $ \Pic(C,\ZZ_{\rm fin}) \to \Iso(\ZZ)$ is surjective. 

Suppose $\a_F = 1$, i.e., 
$F$ is an auto-equivalence of $\Gr(C,\ZZ_{\rm fin})$ such that $FX_n \oplus FY_n \cong X_n 
\oplus Y_n$ for all $n \in \ZZ$.  

Suppose $FC \cong C(I)$. Since $\hom(C,Y_n)=0$ for all $n$,
$\hom(C(I),FY_n)=0$ for all $n \in \ZZ$. However, $\hom(C(I),Y_n) \ne 0$ if $n \in I$ so, if
$n \in I$, then $FY_n \cong X_n$. If $n \notin I$, then $\hom(C(I),X_n) \ne 0$ so $FY_n\cong Y_n$ if
$n \notin I$. Hence $FY_n \cong X_n$ if and only if $n \in I$. 
But $Y_n(I) \cong X_n$ if and only if $n \in I$ so $FY_n \cong Y_n(I)$ for all $n \in \ZZ$.
It follows that $FX_n \cong X_n(I)$ for all $n \in \ZZ$. Since $FS \cong S(I)$ for all special simples $S$, $F \cong (I)$. 
\end{pf}

There is a $\ZZ$-linear action of $\Pic(\gr C)$ on $K_0 (C,\ZZ_{\rm fin})$ given by
$$
[F]\cdot[M]:=[FM].
$$
It is simpler to write this as  $F.[M]:= [FM]$.  

Let $\fp$ denote the  kernel of the  rank function $$\rank:K_0(\gr C) \to \ZZ.$$
Since auto-equivalences preserve rank,  $\fp$ is stable under the action of $\Pic(C,\ZZ_{\rm fin})$. 
Because the rank function is surjective $\fp$ is a prime ideal and 
$$
K_0(\gr C) = \fp \oplus \ZZ\cdot[C].
$$
 
 \begin{prop}
 \hfill
\begin{enumerate}
  \item 
The elements $[X_n]$, $n \in \ZZ$, are a basis for $\fp$.
   \item 
   The elements $[X_n]$ and $[Y_n]$, $n \in \ZZ$ form a full set of pairwise distinct 
   representatives of $\fp/2\fp$.
   \item{}
  $\big\{[X_n],[Y_n] \; \big| \;n \in \ZZ\big\}$ and $\fp/2\fp$ are $\Pic(C,\ZZ_{\rm fin})$-torsors.
     \item 
 $\fp^2=\fp$.
   \end{enumerate}
 \end{prop}
 \begin{pf}
 (1)
If we identify $\ZZ\ZZ_{\rm fin}/\ker \Upsilon$ with its image in $K_0(\gr C)$, then 
$$
1-u_n = [X_n] = - [Y_n]
$$
for all $n \in \ZZ$.  Since $\{1,u_n \; | \; n \in \ZZ\}$ is a basis for $\ZZ\ZZ_{\rm fin}/\ker \Upsilon$
and $\rank(X_n)=0$,  the elements $[X_n]$, $n \in \ZZ$,   form a basis for $\fp$.
 
(2)
This follows immediately from (1). 

(3)
By Theorem \ref{thm.symm}, $\Pic(C,\ZZ_{\rm fin})$ acts transitively on the set
 $\{[X_n],[Y_n] \; | \;n \in \ZZ\}$. An auto-equivalence is determined up to isomorphism by its action
on the special simples (Corollary \ref{cor.FG}) so the only auto-equivalence that acts trivially on 
$\{[X_n],[Y_n] \; | \;n \in \ZZ\}$ is the identity.

 (4)
 This follows from the fact that  $X_i \otimes Y_i \cong Y_i$ and
 $Y_i \otimes Y_i \cong X_i$. 
 \end{pf}

\section{The correspondence between $\Gr A$ and $\Gr(C,\ZZ_{\rm fin})$}
\label{sect.A+C}

In this section we examine the correspondence implemented by the equivalence 
$\Hom(P,-)$  between various significant features of $C$  and $A$. 
The key to doing this this is to match up the special simple $C$-modules 
with the corresponding simple $A$-modules.

 \subsection{The special simple graded $A$-modules}
 Following Sierra we define  the graded simple $A$-modules 
 $$
 X:= {{A}\over{xA}} \qquad Y:=\biggl({{A}\over{yA}}\biggr)(-1).
 $$
 We call the $X(n)$s and $Y(n)$s, $n \in \ZZ$, the {\sf special} simple $A$-modules. 
 
 Recall that $A_0=k[xy]$.

\begin{prop}
\cite[Lemma 4.1]{Sue1}
\label{prop.ss4.1}
The simple graded $A$-modules are
\begin{enumerate}
  \item 
  the modules $X(n)$ and $Y(n)$, $n \in \ZZ$, and 
  \item 
 the modules $A/\fm A$ where $\fm$ is a maximal ideal of $A_0$ but not  one of the ideals $(xy-n)A_0$
  for any  $n \in \ZZ$.
\end{enumerate}
\end{prop}

 We note that 
 $$
 X(n)_m  \ne  0 \Longleftrightarrow m \le -n
 $$
 and 
  $$
 Y(n)_m  \ne  0 \Longleftrightarrow m \ge -n+1
 $$
 whereas the non-special simple graded $A$-modules are non-zero in all degrees.
 It follows that the special simple $A$-modules can be recognized by the degrees in which they are non-zero.

 \begin{prop}
 \label{prop.match.simples}
 Let $H(P,-):  \Gr(C,\ZZ_{\rm fin}) \to \Gr A$  be the equivalence  in the proof of Theorem \ref{thm.main}. Then
 $$
 H(P,X_n) \cong 
 \begin{cases}
 Y(n) & \text{if $n>0$,}
 \\
 X(n) & \text{if $n \le 0$,}
 \end{cases}
 $$
 and 
  $$
 H(P,Y_n) \cong 
 \begin{cases}
 X(n) & \text{if $n>0$,}
 \\
 Y(n) & \text{if $n \le 0.$}
 \end{cases}
 $$
 \end{prop}
 \begin{pf}
 Since $X_n$ is simple, so is $H(P,X_n)$.
 Now
 $$
 H(P,X_n)_m = \hom(P_{(-m,*)},X_n) \cong \hom(C ([-m]),X_n) \cong (X_n)_{[-m]}.
 $$
 Therefore  
 $$
  \hom(P, X_n)_m \ne 0 \Longleftrightarrow n \notin [-m]   \Longleftrightarrow 
  \begin{cases}
  		m \ge -n+1  & \text{if $n > 0$}
		\\
		 m \le -n & \text{if $n \le 0$.}
	\end{cases}
$$
 Hence  $H(P,X_n)$ is as described.  The argument for $H(P,Y_n)$ is similar.  
 \end{pf}

  Let $\s$ be the automorphism of $A$  defined by $\s(x)=y$ and $\s(y)=-x$. 
  For every $n \in \ZZ$,  $\s_*(X(n))$ is isomorphic to $A/yA$ as an ungraded $A$-module and 
   $\s_*(Y(n))$ is isomorphic to $A/xA$ as an ungraded $A$-module.
 But $\sigmaol(m) = -m$ for all $m \in \ZZ$, so 
 \begin{equation}
 \label{eq.sigma*.simples}
  \s_*(X(n))    \cong Y({-n+1})
  \qquad \hbox{and}
  \qquad 
  \s_*(Y(n))\cong X({-n+1})
 \end{equation}
 for all $n \in \ZZ$. 
   
\begin{prop}
\label{prop.fourier}
Let $P$ be the bigraded $A$-$C$-bimodule in the proof of Theorem \ref{thm.main}.
Then 
$$
\sigma_*  \circ H(P,-) \cong  H(P,-) \circ \tau_*\varphi_*.
$$
\end{prop}
\begin{pf}
By Corollary \ref{cor.FG}, it suffices to show that  
$\sigma_* H(P,S) \cong H(P,\tau_*\varphi_*S)$ for every special simple $S$. 
 
If $n >0$, then  $\sigma_* H(P,X_n) \cong   \s_*(Y(n)) \cong X(-n+1)$ and 
 $\sigma_* H(P,Y_n) \cong   \s_*(X(n)) \cong Y(-n+1)$. If $n \le 0$, 
 then  $\sigma_* H(P,X_n) \cong   \s_*(X(n)) \cong Y(-n+1)$ and 
 $\sigma_* H(P,Y_n) \cong   \s_*(Y(n)) \cong X(-n+1)$.
 
 Now we consider the action  of $H(P,-) \circ \tau_*\varphi_*$. 
 By Theorem \ref{thm.tau*.simples}, $\tau_*\varphi_* X_n \cong X_{-n+1}$ and 
 $\tau_*\varphi_* Y_n \cong Y_{-n+1}$. If $n \le 0$, then $-n+1 >0$ so, by
 Proposition \ref{prop.match.simples},
 $H(P,-) \circ \tau_*\varphi_*X_n \cong Y(-n+1)$ and  $H(P,-) \circ \tau_*\varphi_*Y_n \cong X(-n+1)$.
 If $n>0$, then $-n+1 \le 0$ so , by
 Proposition \ref{prop.match.simples},
 $H(P,-) \circ \tau_*\varphi_*X_n \cong X(-n+1)$ and  $H(P,-) \circ \tau_*\varphi_*Y_n \cong Y(-n+1)$.
Comparing the results in this and the previous paragraphs, it follows that $\sigma_*  \circ H(P,-) \cong  H(P,-) \circ \tau_*\varphi_*$.
 \end{pf}

\begin{lem}
Let
$$
\Sigma := (\{1\}) \circ \tau_* .
$$
Then $\Sigma$ is an automorphism of $\Gr(C,\ZZ_{\rm fin})$ and it permutes the isomorphism classes of the special simple modules as in the diagram:
$$
\UseComputerModernTips
  \xymatrix{
  \cdots \ar[r] & X_{-1}  \ar[r] & X_{0} \ar[dr] & X_{1}  \ar[r] & X_{2}  \ar[r] & \cdots 
  \\
 \cdots  \ar[r] & Y_{-1}  \ar[r] & Y_{0} \ar[ur] & Y_{1}  \ar[r] & Y_{2}  \ar[r] & \cdots 
}
$$
\end{lem}
\begin{pf}
 By Theorem \ref{thm.tau*.simples}, $\Sigma X_n \cong X_{n+1}(\{1\})$. Hence $\Sigma X_0 \cong Y_1$.
 By the remarks at the beginning of section  \ref{sect.specials}, if $n \ne 0$, then  $\Sigma X_n \cong X_{n+1}$.
 Similar considerations apply to $\Sigma Y_n$. 
\end{pf}

\begin{prop}
\label{prop.shift.Sigma}
Let $H(P,-):\Gr(C,\ZZ_{\rm fin}) \to \Gr A$ be the equivalence in Theorem \ref{thm.main}. 
Then there is an isomorphism  of functors $$H(P,-) \circ \Sigma \cong (1) \circ H(P,-).$$
\end{prop}
\begin{pf}
By Corollary \ref{cor.FG}, it suffices to show that  
$H(P,S)(1) \cong H(P,\Sigma S)$ for every special simple $S$. 
 \end{pf}

 \subsection{The automorphisms $\iota_J$}
 \label{sect.iota}

 The key to much of Sierra's analysis of $\Gr A$ is her discovery of the automorphisms  
 $\iota_J$,  $J \in \ZZ_{\rm fin}$, of $\Gr A$ that she describes in  \cite[Prop. 5.9]{Sue1}. 
 These have the following properties: 
\begin{enumerate}
  \item 
  $\iota_\varnothing=\id_{\Gr A}$;
  \item 
  $\iota_n = (n) \circ \iota_0 \circ(-n)$ for all $n \in \ZZ$;\footnote{Because our $(n)$ is equal to Sierra's $\langle -n\rangle$, our $\iota_n$ is her $\iota_{-n}$.}
  \item 
  $\iota_J = \prod_{j \in J} \iota_j$ for all $J \in \ZZ_{\rm fin}$;
  \item
  $\iota_0^2 \cong \id_{\Gr A}$. This is not an equality.
\end{enumerate}
The third condition says that the map $J \to \iota_J$
 from $\ZZ_{\rm fin}$ to the automorphism group of $\Gr A$
is a homomorphism of monoids, and the fourth condition implies that 
the map $\ZZ_{\rm fin} \to \Pic(\Gr A)$, $J \mapsto \iota_J$, or, more precisely, $J \mapsto$ the image of $\iota_J$ in $\Pic(\Gr A)$,  is a group homomorphism.

The functor $\iota_0$ is defined first as an automorphism of the subcategory of $\Gr A$ consisting of the 
projective graded modules and as an automorphism  of that subcategory it is a subfunctor of the identity functor. It follows that every $\iota_J$ is also a  subfunctor of the identity functor on that subcategory.

 Sierra shows that $\hom(P,X \oplus Y) \cong k$ for all rank one graded projectives $P$ (cf. Corollary \ref{cor.proj=free} and parts (5) and   (6) of Proposition  \ref{prop.simples}). The functor $\iota_0$ is
 then defined on a rank one projective by 
 \begin{align*}
 \iota_0P :=& \; \ker f \;  \; \hbox{ where $f:P \to X \oplus Y$ is any}
 \\
 & \qquad \quad \;  \hbox{non-zero graded homomorphism}.
 \end{align*}
Equivalently, $\iota_0P$ is the unique graded submodule of $P$ that fits into an exact sequence
$$
0 \to \iota_0 P \to P \to X \oplus Y
$$
in which the right-most map is the unique (up to scalar multiple) non-zero map $P \to X \oplus Y$.
From the exact sequence $0 \to \iota_0(P(-n)) \to P(-n) \to X \oplus Y$, we see that
$\iota_nP$ is the unique submodule of $P$ fitting into an exact sequence
$$
0 \to \iota_n P \to P \to X(n) \oplus Y(n)
$$
where the right-most map is non-zero.

\begin{thm}
\label{thm.iotaJ.CJ}
Let $J \in \ZZ_{\rm fin}$.
Then 
$$
H(P,-) \circ (J) \cong \iota_J \circ H(P,-).
$$
\end{thm}
\begin{pf}
By \cite[Prop. 5.9]{Sue1}, $\iota_n$ interchanges the isomorphism classes of $X(n)$ and $Y(n)$ and 
fixes the isomorphism classes of all other $X(m)$s and $Y(m)$s.  On $\Gr(C,\ZZ_{\rm fin})$,
the twist $(\{n\})$ interchanges $X_n$ and $Y_n$ and fixes the isomorphism classes of all other special
simples. It now follows from Proposition \ref{prop.match.simples} that   $ \iota_n H(P,S)
\cong H(P,S(\{n\}))$   for all special $C$-modules $S$. 
Hence by Corollary \ref{cor.FG}, $ \iota_n \circ H(P,-) \cong H(P,-) \circ (\{n\}) $. 
Thus, if $J=\{j_1,\ldots,j_t\}$, then
\begin{align*}
 \iota_J \circ H(P,-) & \cong  \iota_{j_1}\circ \cdots  \circ \iota_{j_t} \circ H(P,-) 
 \\
 & \cong H(P,-) \circ (\{j_1\}) \circ  \cdots  \circ  (\{j_1\}) 
 \\
 &  \cong  H(P,-) \circ (J),
 \end{align*}
 as required. 
\end{pf}

\subsection{The monoidal structures}
The equivalence of categories does not respect the ``natural'' internal tensor products on $\Gr A$ and 
$\Gr(C,\ZZ_{\rm fin})$.  
The tensor product on $\Gr(C,\ZZ_{\rm fin})$ is the usual graded tensor product over $C$. 
 The  tensor product on $\Gr A$ is that which exists on the category of $\cD_Y$-modules  for any 
 smooth variety $Y$, namely $\cM \otimes \cN= \cM \otimes_{\cO_Y} \cN$ with a derivation 
 $\d$ acting on the tensor product as  $\d \otimes 1 + 1 \otimes \d$. 
 Specializing to $A$ and taking $k[y]$   as the coordinate ring  of the line on which $A$ acts as differential operators (and $x$ acts as  $-d/dy$), the internal tensor product on $\Gr A$ is $-\otimes_{k[y]} -$ with $x$ 
 acting on the tensor product as $x \otimes 1 + 1 \otimes x$. 

Both tensor products are commutative and  
\begin{align*}
X(m) \otimes  X(n) & \cong X(m+n),
\\
Y(m) \otimes_{} Y(n) & \cong Y(m+n-1)
 \cong  X(m-1) \otimes  Y(n) 
 \end{align*}
whereas
$$
X_i \otimes X_j \cong Y_i \otimes Y_j \cong \d_{ij} X_i, \quad X_i \otimes Y_j \cong \d_{ij} Y_j.
$$
The identity for the tensor product in $\Gr A$ is the simple module $X=A/xA$,
 whereas the identity for the tensor product in $\Gr(C,\ZZ_{\rm fin})$ is the projective module $C$.

The tensor product of two finitely generated $C$-modules is finitely generated but that property 
does not hold for $A$-modules.

 \subsection{Preparations for section \ref{sect.tw.hcr}}
 \label{sect.preps}
 
 \subsubsection{The left $k[z]$ action on graded right $A$-modules}
 
 Throughout we will write 
 $$
 z=xy.
 $$
 The degree zero component of $A$ is therefore
 $$
 A_0=k[z].
 $$
 
 Sierra \cite[Sect. 4, p.14]{Sue1} makes the fundamental observation that 
 every graded right $A$-module $M$ can be given the structure 
 of a {\it left} $k[z]$-module in such a way that $M$ becomes a $k[z]$-$A$-bimodule. 
 The left action of $z$ on an element $m \in M_j$ is 
 \begin{equation}
 \label{rt.z.action}
z .  m:=m(z-j).
 \end{equation}
 The left $k[z]$-action commutes with the right $A$ action because for all homogeneous
 $a \in A$,
 \begin{equation}
 \label{conj.z}
 [a,z]=   (\deg a) a. 
 \end{equation}
Homomorphisms  $f:M \to N$ in $\Gr A$ are also homomorphisms of left $k[z]$-modules.
Hence $\hom(M,N)$ has induced left and right $k[z]$-module structures given by 
$$
(z.\psi)(m)= z.\psi(m) \quad \hbox{and} \quad  
(\psi.z)(m)=\psi(z.m).
$$
Since $\psi$ preserves degree these two $k[z]$-module structures are the same.

When $M=A$ the left action of $z$ on $A$ given by (\ref{rt.z.action}) is the ordinary 
left multiplicative action, i.e., $z.a=za$.  This follows from (\ref{conj.z}).
However,  the left action of $z$ on $A(n)$ given by (\ref{rt.z.action}) coincides with left
multiplication by $z+n$.  For example, if $\1ol$ denotes 1 viewed as an element in $A(n)$ then
$z.\1ol =\1ol (z+ n)$ because $\1ol \in A(n)_{-n}$.

If $z.M=0$, then $(z-n).M(n)=0$. 

Recall that $\cO_\l=A/(z-\l) A$. It is straightforward to see that 
$$
(z-n).\cO_n= (z-n).\cO_{n+1}(-1)= 0
$$
and, as a consequence of either of these facts,
 $$
 (z-n) .X(n)= (z-n).Y(n)=0.
 $$

\subsubsection{Isomorphisms between products of the $\iota_J$s}
\label{sect.eta}

 For each $J \in \ZZ_{\rm fin}$, define
 \begin{equation}
 \label{defn.fJ}
 h_J=\prod_{j\in J}(z-j).
 \end{equation}
This polynomial belongs to $k[z]=A_0$.

Because the left action of $(z-n)$ annihilates  the non-split extensions between $X(n)$ and $Y(n)$
 it follows that 
$$
\iota_n^2P=(z-n)P.
$$
There is therefore a unique isomorphism 
$$
\eta^n:\iota_n^2 \to \id_{\Gr A}
$$ 
such that $(\eta^n)_P: \iota_n^2P \to P$ is left multiplication  by $(z-n)^{-1}$ for every projective $P$.
To be precise, $\eta^n$ is first defined as a natural transformation between the restrictions of the functors to the subcategory of projectives, and as such $\eta^n$ is multiplication by 
$(z-n)^{-1}$. The natural transformation $\eta^n$ ``extends'' uniquely to a natural transformation between the functors defined on all of $\Gr A$ {\it but}, as a natural transformation on $\Gr A$,
 $\eta^n$ is {\it not} ``multiplication by $(z-n)^{-1}$.''  Similarly, if $I,J \in \ZZ_{\rm fin}$, we define  the isomorphism 
$$
\eta^{IJ} : \iota_{I \oplus J}   \to \iota_I\iota_J 
$$
to be left multiplication by the polynomial $h_J$ in (\ref{defn.fJ}).

 We define the automorphism $\s_n:\iota^2_n A \to A$ to be the isomorphism $\eta^n$ at $A$, i.e.,
 $\s_n=(\eta^n)_A$. Thus $\s_n$ is left multiplication by $(z-n)^{-1}$. We also define
 \begin{equation}
 \label{sigmaJ}
 \s_J:=\prod_{j \in J} \s_j \; : \; \iota_J^2 A \to A
 \end{equation}
 and $\s_\varnothing=\id_A$.

 \section{$C$ is a twisted homogeneous coordinate ring for $\Gr A$}
 \label{sect.tw.hcr}
 
 In this section $C$ is constructed directly from $A$ as a sort of twisted homogeneous coordinate
 ring for $\Gr A$. We will show that $C$ is isomorphic as a graded ring to
  the ring $B$ defined in (\ref{defn.B}).

 \subsection{The $\ZZ_{\rm fin}$-graded ring $B$}

 We define the $\ZZ_{\rm fin}$-graded
 ring 
\begin{equation}
\label{defn.B}
B := \bigoplus_{J \in \ZZ_{\rm fin}}  B_J = \bigoplus_{J \in \ZZ_{\rm fin}}  \hom(A,\iota_J A)
\end{equation}
endowed with the following multiplication: if $f \in B_I$ and $g \in B_J$, then
\begin{equation}
\label{eq.B.prod}
f\cdot g:= \iota_{I \oplus J}(\s_{I \cap J}) \circ  \iota_J(f) \circ g.
\end{equation}
where $\sigma_{I \cap J}$ is defined in (\ref{sigmaJ}).

 The identity component of $B$ is 
 $$
 B_\varnothing=\hom(A,A) \cong A_0 = k[z].
 $$

  \begin{lem}
  \label{lem.assoc.B}
The product (\ref{eq.B.prod}) on $B$ is associative.
\end{lem}
\begin{pf}
It suffices to check that the natural transformations $\eta^{IJ}$ defined in section \ref{sect.eta}
 satisfy the conditions 
mentioned after \cite[Prop. 2.2]{Sue1}. This reduces to checking that the analogue of the commutative 
diagram \cite[(2.5)]{Sue1} really does commute, and that reduces to showing that
$$
\iota_I\bigl( \eta^{J,L}\bigr) \circ \eta^{I,J \oplus L} 
= \eta^{IJ}\iota_L \circ \eta^{I \oplus J,L}
$$
which, in turn, reduces to showing that 
$$
h_{J \cap L} h_{I \cap (J \oplus L)} = h_{I \cap J} h_{(I \oplus J) \cap L}.
$$
This equality follows from the fact that 
$$
({J \cap L}) \cap{I \cap (J \oplus L)} =   ({I \cap J}) \cap {(I \oplus J) \cap L} = \varnothing
$$
and 
$$ 
({J \cap L}) \cup ({I \cap (J \oplus L)}) =   ({I \cap J}) \cup ({(I \oplus J) \cap L} ).
$$
The last two expressions both equal $(I \cap J) \cup (I \cap L)  \cup (J\cap L)$.
\end{pf}

  \begin{lem}
  \label{lem.key}
  With the above notation,
   \begin{enumerate}
  \item 
   $ \bigl( \iota_J A \bigr)_0 =  h_J k[z]$ where $h_J$ is as defined in (\ref{defn.fJ}).
  \item 
  The map 
  $$
  \rho_J:(\iota_JA)_0 \to \hom(A,\iota_JA), \quad \rho_J(m)(a)=ma,
  $$
  is an  isomorphism of $A_0$-modules. 
  \item 
    As a right $B_\varnothing$-module, $B_J$ is freely generated by the elements
    $$
    b_J:=\rho_J(h_J).
    $$
    \item{}
     $ b_n^2 +n =b_m^2+m$  for all $m,n \in \ZZ$. 
     \item
     $b_Ib_J=b^2_{I \cap J} b_{I \oplus J}$.
     \item{}
     $b_J=\prod_{j \in J} b_j$.
     \item{}
     $B$ is a commutative $k$-algebra generated by $\{b_n \; | \; n \in \ZZ\}$.
\end{enumerate}
    \end{lem}
  \begin{pf}
  (1)
   As noted in the proof of \cite[Lemma 5.14]{Sue1},  
\begin{equation}
\label{iotanA}
 \iota_n A = 
 \begin{cases}
 x^{n+1}A + (xy+n)A & \text{if $n \ge 1$}
 \\
  xA & \text{if $n = 0$}
 \\
  yA & \text{if $n = -1$}
 \\
  y^{-n}A + (xy+n)A & \text{if $n \le -2 $}.
 \end{cases}
 \end{equation}
 Hence $ \bigl( \iota_n A \bigr)_0 =  (z-n)k[z]$ and the result follows from the fact that
 $$
 \iota_J A= \bigcap_{j \in J} \iota_j A.
 $$

(2)
  This is trivial.
 We will use the isomorphism $\rho_\varnothing: A_0 = k[z] \to B_\varnothing $ to identify $k[z]$ with $B_\varnothing $.

 (3)
The multiplication $B_J \times B_\varnothing\to B_J$ sends $(f,g)$ to $f \circ g$. Since $B_\varnothing =k[z]$ and 
$B_J=\rho_J(h_J k[z])$ the result follows.  
 
 (4)
The multiplication $B_n \times B_n \to B_\varnothing $  is given by
 $$
 b.b'=\sigma_n \circ  \iota_n(b) \circ b'
 $$
so
 $$
 (b_n.b_n)(a)= \sigma_n(h_n^2a)=(z-n)^{-1}h^2_n a = h_n a \in B_\varnothing .
 $$
Hence $b_n^2 = \rho_\varnothing (h_n) =\rho_\varnothing (z-n) = \rho_\varnothing(z)-n$ and  
 $$
 b_n^2 +n = \rho_\varnothing(z)= b_m^2+m
 $$
 for all $m,n \in \ZZ$. 
 
 (5) 
 By definition of the product in $B$, 
 \begin{align*}
 b_Ib_J & = \iota_{I \oplus J}(\s_{I \cap J})  \circ  \iota_J(b_I) \circ b_J
 \\
 & = \iota_{I \oplus J}(\s_{I \cap J})  \circ  \iota_J(\rho_I(h_I)) \circ \rho_J(h_J)
 \end{align*}
 But $\rho_J(h_J)$ is ``left multiplication by $h_J$'', $\rho_I(h_I)$ is  ``left multiplication by $h_I$'',  $ \iota_J(\rho_I(h_I))$ is the restriction of  $\rho_I(h_I)$  so is also  ``left multiplication by $h_I$'', and 
 $\s_{I \cap J}$, and hence $ \iota_{I \oplus J}(\s_{I \cap J})$, is ``left multiplication by $h^{-1}_{I \cap J}$''.
 Hence 
 \begin{align*}
 b_I b_J  & = \bigl( \hbox{left multiplication by $h_{I \cup J}$}\bigr)  :  A \to \iota_{I \oplus J} A
 \\
 & = \rho_{I \oplus J} (h_{I \cup J}).
 \end{align*}
Hence  $b_{I \cap J} b_{I \oplus J} = \rho_{I \cup J} (h_{I \cup J})=b_{I \cup J}$
  and therefore
   \begin{align*}
 b^2_{I \cap J} b_{I \oplus J}  & =   b_{I \cap J} b_{I \cup J}
 \\
 & = \rho_{(I \cap J) \oplus (I \cup J)}  (h_{(I \cap J) \cup(I \cup J)})
 \\
  & = \rho_{I \oplus J}(h_{I \cup J})
 \\
 & = b_I b_J
 \end{align*}
 as claimed.
 
 (6) 
 This follows from (5) and an induction argument on $|J|$.
 
 (7)
 It follows from (5) that $b_mb_n=b_nb_m$ for all $m, n \in \ZZ$. 
 \end{pf}

 \begin{prop}
 As $\ZZ_{\rm fin}$-graded $k$-algebras,  $B \cong C$.
 \end{prop}
 \begin{pf}
 By parts (5) and (6) of Lemma \ref{lem.key}, the function $\psi(x_n):=b_n$ extends to a well-defined 
 homomorphism  $\psi:C \to B$ of $\ZZ_{\rm fin}$-graded $k$-algebras. 
 By part (7) of Lemma \ref{lem.key}, $\psi$ is surjective. 
 
 Since $C_\varnothing $ is a polynomial ring in one variable  the restriction of $\psi$
 to $C_\varnothing $ is an isomorphism $C_\varnothing  \to B_\varnothing$. 
 The kernel of $\psi$ is the sum of its homogeneous components. If $C_J \cap \ker \psi$ were non-zero
 multiplying it by $x_J$ would produce a non-zero element of $C_\varnothing  \cap \ker\psi$; but the latter is zero, so $C_J \cap \ker \psi =0$ and we conclude that  $\psi$ is an isomorphism, as claimed. 
 \end{pf}

 \begin{lemma}
 [Sierra]
  \cite[Prop. 4.1]{Sue1}\footnote{This result is also a consequence of the equivalence $\Gr A \equiv \Gr (C,\ZZ_{\rm fin})$ and the fact that the $C(J)$s are a set of generators for $\Gr(C,\ZZ_{\rm fin})$.}
 The set of modules $\{\iota_J A \; |\; J \in \ZZ_{\rm fin}\}$ generates $\Gr A$.
 \end{lemma}

 Because the isomorphisms $\eta^{IJ}:\iota_{I \oplus J} \to \iota_I \iota_J$ satisfy the conditions 
 verified in the proof of Lemma  \ref{lem.assoc.B} there is a well-defined functor
 \begin{align*}
F: &\; \Gr A \to \Gr(B,\ZZ_{\rm fin}), 
\\
FM:= & \bigoplus_{J \in \ZZ_{\rm fin}}        \hom(A,\iota_JM).
\end{align*} 
Because $\{\iota_J A \; | \; J \in \ZZ_{\rm fin}\}$ is a set of projective generators for $\Gr A$,
it follows from del Rio's result (Theorem \ref{thm.dR}) that $F$ is an equivalence.

 \subsection{Final remarks}
Is there an a priori reason why the Weyl algebra with the given $\ZZ$-grading 
might be so intimately related to a ring like 
$C$ (or, equivalently, a stack like $\cX$)?
One explanation is this. A $\ZZ$-grading typically forces graded modules 
to behave somewhat like ungraded modules over a ring of dimension one less.  Since the Weyl algebra has Gelfand-Kirillov dimension two, and since rings of Gelfand-Kirillov dimension one behave a lot like curves,
$\Gr A$ might reasonably be expected to  exhibit curve-like features. The stacky behavior corresponds to the 
existence of non-split extensions between non-isomorphic simples.

  \end{document}